\documentclass[reqno,11pt]{amsart}
\usepackage{amsmath, amssymb, amsthm,amsfonts} 
\usepackage[english]{babel}
\usepackage{bbm,bm}
\usepackage{graphicx}
\usepackage{url}
\usepackage{epstopdf}
\usepackage[a4paper,bindingoffset=0.2cm,left=1cm,right=1cm,top=2.5cm,bottom=2cm,footskip=.8cm]{geometry}
\usepackage{subcaption}
\usepackage{textcomp}
\usepackage[]{algorithm2e}
\usepackage{multicol}
\usepackage{multirow}
\usepackage{wrapfig}
\usepackage{diagbox}
\usepackage{floatrow}
\newfloatcommand{capbtabbox}{table}[][\FBwidth]

\usepackage{tikz}
\usetikzlibrary{arrows, automata,positioning,calc,shapes,decorations.pathreplacing,decorations.markings,shapes.misc,petri,topaths}
  \usepackage{pgfplots}
  \pgfplotsset{compat=newest}
  \usetikzlibrary{plotmarks}
  \usepackage{grffile}
\newlength\figureheight
  \newlength\figurewidth
\pgfplotsset{%
    tick label style={font=\scriptsize},
    label style={font=\footnotesize},
    legend style={font=\footnotesize},
         every axis plot/.append style={very thick}
}

\usepackage{rotating}
\usepackage{amsbsy,enumerate}
\usepackage{graphicx}
\usepackage{comment}
\usepackage{mathrsfs} 
\usepackage{tikzscale}
\usepackage{xcolor}

\newcommand{\bs}{\boldsymbol}

\newcommand{\vb}{\vspace{5mm}}

\newcommand{\minus}{\scalebox{0.5}[1.0]{$-$}}

\DeclareMathOperator*{\argmin}{arg\,min}

\DeclareMathOperator{\NB}{NB}
\DeclareMathOperator{\Exp}{Exp}

\allowdisplaybreaks

\makeatletter
\newtheorem*{rep@theorem}{\rep@title}
\newcommand{\newreptheorem}[2]{%
\newenvironment{rep#1}[1]{%
 \def\rep@title{#2 \ref{##1}}%
 \begin{rep@theorem}}%
 {\end{rep@theorem}}}
\makeatother

\makeatletter
\renewcommand*\env@matrix[1][\arraystretch]{%
  \edef\arraystretch{#1}%
  \hskip -\arraycolsep
  \let\@ifnextchar\new@ifnextchar
  \array{*\c@MaxMatrixCols c}}
\makeatother

\newtheorem{theorem}{Theorem}
\newreptheorem{theorem}{Theorem}

\newtheorem{corollary}{Corollary}
\newtheorem{assumption}{Assumption}
\newtheorem{remark}{Remark}
\newtheorem*{example}{Example}

\newtheorem{proposition}{Proposition}
\newreptheorem{proposition}{Proposition}

\newcommand{\edsger}{\mbox{\sc edsger}^\star }
\newcommand{\edsgernonstar}{\mbox{\sc edsger} }

\begin{document}


\title[A Framework for Efficient Dynamic Routing 
under Stochastically Varying Conditions]
{A Framework for Efficient Dynamic Routing 
\\ under Stochastically Varying Conditions}

\author{Nikki Levering, Marko Boon, Michel Mandjes, and Rudesindo N\'u\~nez-Queija}

\maketitle




%



\begin{abstract}
Despite measures to reduce congestion, occurrences of both recurrent and non-recurrent congestion
cause large delays in road networks with important economic implications.
Educated use of Intelligent Transportation Systems (ITS) can significantly reduce travel times.
We focus on a dynamic stochastic shortest path problem: our objective is to minimize the expected travel time of a
vehicle, assuming the vehicle may adapt
the chosen route while driving. We introduce a new stochastic process that incorporates ITS information to model
the uncertainties affecting congestion in road networks. A Markov-modulated background process tracks
traffic events that affect the speed of travelers. The resulting continuous-time routing model
allows for correlation between velocities on the arcs and incorporates both recurrent and non-recurrent congestion.
Obtaining the optimal routing policy in the resulting semi-Markov decision process 
using dynamic programming
is computationally intractable for realistic network sizes. To overcome
this, we present the $\edsger$ algorithm, a Dijkstra-like shortest path algorithm that can be used dynamically with
real-time response. We develop additional speed-up techniques that reduce the size of the network model. We quantify the performance of the algorithms by providing numerical examples that use
road network detector data for The Netherlands.

\vb

\noindent
{\sc Keywords.} Routing $\circ$ real time $\circ$ shortest path $\circ$ Semi-Markov decision processes $\circ$ dynamic programming.

\vb

\noindent
{\sc Affiliations.} 
{NL, MM, and RNQ  are with the Korteweg-de Vries Institute for Mathematics, University of Amsterdam, Amsterdam. MB is with the Department of Mathematics and Computer Science, Eindhoven University of Technology, Eindhoven. MB and MM are also with E{\sc urandom}, Eindhoven. MM is also with Amsterdam Business School, University of Amsterdam, Amsterdam. RNQ is also with Centrum Wiskunde \& Informatica, Amsterdam.}
This research project is partly funded by the NWO Gravitation project N{\sc etworks}, grant number 024.002.003. Date: \today. 

\vb

\noindent
{\sc Declaration of conflicts of interest:} none.

\end{abstract}

\newpage

\section{Introduction}

To reduce congestion in road networks, 
a wide variety of measures can be thought of. 
Perhaps the most basic measure is to increase the network's capacity --- the Dutch government 
for example increased the total road length from 130.446 km in 2001 to 141.361 km  in 2020 \cite{CBS}.
Alternatively, one may attempt to better exploit the existing resources, examples being
the use of reversible lanes \cite{Wolshon2006} and the
transit lane experiments in various US states \cite{Boriboonsomsin2008, Chang2008}.
Despite such measures, 
recurrent congestion (i.e., congestion during peak hour)
and non-recurrent congestion (i.e., congestion due to incidents) 
remain a major concern. Consequences include
substantial
delay in travel times, 
increased economic costs, and 
negative environmental effects.

So as to reduce the effects of recurrent congestion,
historical data is used to incorporate periodically occurring 
events (rush hours, weekly patterns, etc.)~into routing algorithms.
Regarding non-recurrent congestion, an important role is played by
Intelligent Traffic Systems
(ITS)
capable of providing 
travelers
with real-time information.
The availability of data on virtually all thinkable network characteristics, 
combined with instant communication to individual travelers, 
offers the potential for routing with minimal delay and minimal congestion.
A complication, however, lies in the unpredictability of traffic surges and incidents,
and in addition in the computational 
effort needed to process all available information in real time.
Our work provides an approach that handles both these challenges: 
providing fast and close-to-optimal routing, using a probabilistic framework
that is well suited for modeling uncertainties that affect congestion.

{{\sc Routing framework}} \\
In this paper we study routing by casting it as a dynamic stochastic 
shortest path problem. Our objective is to minimize
the expected travel time of a vehicle pertaining to a given 
origin-destination pair, 
assuming that the vehicle is allowed to adapt the chosen route while driving.
A key feature of our setup is that we model the 
travel-time dynamics using a suitably chosen stochastic processes.
Specifically, for each road of the network, we let 
the velocity of a vehicle on that road
be determined by the state of a Markovian background process, describing
incidents, weather conditions, etc.

Routing decisions in the proposed dynamic model can 
be expressed in terms of a semi-Markov
decision process (SMDP) \cite{Boucherie2017, Puterman1994, Sutton2018}.
The state of the SMDP consists of the location
of the vehicle, which, in a practical context, can be traced by GPS, 
and the state of the background process, being
provided by an ITS.
The route can be adapted at each intersection along the traveled path.
Based on the state of
the SMDP at an intersection, a decision is taken that informs the driver
which arc to travel next.
Clearly, the use of a Markovian background process facilitates
the modeling of non-recurrent, random events. 
Perhaps surprisingly, however, we will argue that it also allows
us to incorporate
the deterministic patterns corresponding to recurrent, (near-)deterministic events.
Importantly, our framework allows for correlation between the travel times on the edges
and is therefore well capable of modeling the spillback effect:
an incident at an arc
causes a drop in speed levels at upstream arcs.

In the context of our model, an optimal routing policy can be evaluated by the use of
dynamic programming (DP), 
as it can be shown that such a policy satisfies
the Bellman optimality equations. 
DP does however suffer from
the curse of dimensionality, making it prohibitively slow
for networks of practical relevance.
To cope with this, we have developed the {\sc edsger} algorithm, as well as its more efficient variant $\edsger$,
which are Dijkstra-like shortest path algorithms, but with the additional feature that
the route can be adapted along the way, based on the current state of the background process.
We furthermore present speed-up techniques
that greatly reduce the size of the underlying network model and its corresponding state-space.
These techniques can therefore be used as preprocessing steps 
to make the routing algorithms
substantially more efficient.

To make our approach operational, 
the parameters pertaining to the background process must be estimated.
Deterministic patterns and corresponding transitions in the
background process can be identified from historical data. 
In addition, we need to estimate 
the frequency and consequences of
non-recurrent events, such as
the frequency of incidents or the drop in speed
level under bad weather conditions.
Examples of institutions that 
collect data on traffic flows and traffic speeds
are the National Data Warehouse for Traffic Information (NDW) in
The Netherlands, the Mobilitäts Daten Marktplatz (MDM) in Germany; in addition various US states have such an infrastructure (see for example MITS, 
the Michigan Intelligent Transportation Systems).

This paper has the following two main contributions: 
\begin{itemize}
    \item [$\circ$]First, we present
the Markovian background process and corresponding SMDP framework
in detail. To our knowledge, this is the first study that
uses a continuous-time Markovian background process
to develop an adaptive routing algorithm. 
It has important advantages over existing approaches: in a single framework it 
offers the potential of incorporating both recurrent and non-recurrent events and
allows for correlation between speeds on different edges.
In addition, our setup guarantees continuity in travel times (travel times are a continuous function of the departure epoch), implying the so-called FIFO-property (the arrival epoch is an increasing function of the departure epoch).
 \item [$\circ$]
Second, after having presented a way to evaluate the 
optimal policy, which suffers from the curse of dimensionality, we propose efficient alternatives.
Through a series of numerical examples, we show that in particular the algorithm $\edsger$ is fast, in the sense that it can be performed in real time, 
providing suboptimal routes in a very small subset of all possible scenarios.
The examples all concern specific parts of the
Dutch road network, with instances that are based on the
data provided by NDW to make sure they are representative for real traffic scenarios. 
\end{itemize}





{{\sc Literature overview}} \\
Routing algorithms have extensively been examined in literature.
In a deterministic graph with positive arc costs satisfying
the FIFO property, the shortest path can be obtained by
Dijkstra's algorithm, developed by Edsger W. Dijkstra in the
late 1950s \cite{Dijkstra1959, Kaufman1993}.
Throughout the years many variations and speed-up
techniques for Dijkstra's algorithm have been presented,
examples of which are the A$^\star$-algorithm and 
the bidirectional Dijkstra algorithm
\cite{Fu2006, Algorithmics}.
Bellman \cite{bellman1958} and Ford \cite{Ford1956} 
constructed a shortest path algorithm that allows for both 
positive and negative arc costs.

The search for optimal routes
can
complicate considerably
as soon as randomness is introduced. The stochastic shortest path problem is an extension of the
deterministic shortest path problem in which uncertainty
in travel times is taken into account.
Several variations of this
problem exist, each with its own objective.
Examples of two well-studied objectives are the minimization 
of the expected travel-time
and the maximization of the on-time arrival probability.
In the present article we focus on the first objective;
for studies on the latter we refer to 
\cite{Fan2005, Frank1969, Nie2009, Pedersen2019}
and references therein.
In case of a minimum expected travel-time objective,
Dijkstra's algorithm can still be applied if there is neither
correlation (incidents on specific arcs do not influence 
travel times on other arcs) nor time-dependency 
(the travel times do not depend on the hour of the day)
\cite{Loui1983, Mirchandani1986, Murthy1996}.
Recurrent events are well-described by time-dependent velocities (think of predictable, recurring events such as rush hours), but
Hall \cite{Hall1986} showed that 
classic
shortest path algorithms like Dijkstra's algorithm fail
in stochastic time-dependent networks.
Hall additionally argued that it is suboptimal to solely
consider static routing and that
it is advantageous to allow travelers to adapt the chosen route
while traveling.
These findings resulted in studies on optimal adaptive routing
algorithms in stochastic time-dependent networks.
Examples of such algorithms are
presented by Fu and Rilett \cite{Fu1998}
and Miller-Hooks and Mahmassani \cite{Miller-Hooks2000}.
These algorithms do, however, not
take ITS information into account and are therefore
limited in their capability to incorporate non-recurrent congestion.

Algorithms that do take ITS information into
account can be divided into two categories.
The first category consists of algorithms that assume that
ITS provides information on all realized travel times in the
network. Without attempting to provide an exhaustive overview, 
we refer in this context to the algorithms presented in 
\cite{Bander2002,Cheung1998,Davies2003, Gao2006, Gao2012,Polychronopoulos1996}.
The second category consists of algorithms that 
assume that ITS provides information on the state of a
background process in the network.
Psaraftis and Tsitsiklis \cite{Psaraftis1993} were among
the first to work in this setting, 
presenting a framework in which the travel
time on an arc depends on an environmental variable becoming
known once a vehicle arrives at the attached node.
Their setup was extended by Kim et al.\ \cite{Kim2005a},
who work with a Markov process, the state of which
state directly determines the travel time. 
By ITS, the state of a Markov
process of every arc is known while traveling through the network, 
so that we can phrase the optimal routing problem in terms of a  Markov decision process (MDP).
The MDP framework is further explored in
\cite{Guner2012, Sever2013, Thomas2007}, 
under the assumption of independence 
of the Markov processes on the arcs and thus 
neglecting spatial correlation. 
They furthermore assume that the travel
time on a link is known once the 
intersection at the head of the link is reached,
whereas in reality the conditions on an arc may still
change while traveling on that arc. Other drawbacks are (i)~the fact that 
the realized travel times can only take values
from a previously chosen discrete collection of values
and (ii)~the fact that additional assumptions
are needed to guarantee the FIFO-property.
In the MDP framework an optimal policy can be
characterized by the Bellman optimality
equations and can therefore be evaluated by performing a DP recursion
\cite{Bellman1957a, Bellman1957b, bertsekas1976, Bertsekas2005, Puterman1994, Ross1983, Sutton2018}. 
DP suffering from the curse of dimensionality,
it may lead to excessive computational costs.
Possible solutions are dimensionality reduction, 
as described in e.g.\ \cite{Sykora2008} and references therein,
and approximative procedures such as approximate dynamic programming (ADP)
\cite{Bertsekas2012, Powell2007}.
In the context of routing, dimensionality reduction
is employed in \cite{Guner2012, Kim2005b},
while ADP is extensively studied in \cite{Sever2013}.

{{\sc Contributions}} \\
Above, we stated our two main contributions, which we will describe in greater detail now. Regarding the first contribution, 
the specific Markov-modulated
background process that we advocate in this paper, describing the speeds at which vehicles can drive on the arcs, 
overcomes the above-mentioned drawbacks of earlier approaches.
Although the approaches of Ferris and Ruszcy\'nski \cite{Ferris2000} and Karoufeh and Gautam \cite{Karoufeh2004} also work with
continuous-time Markovian processes, there are important differences
with our proposal.
In the setup of \cite{Ferris2000},
the travel times are directly modeled relying on a continuous-time Markov process. 
In \cite{Karoufeh2004}, a continuous-time Markov background process 
is considered to model the
velocities, but it only describes
the dynamics on a single arc and does not consider routing.
With our continuous-time
Markovian background process, applying to the network as a whole, we overcome the drawbacks
of the discrete MDP framework that we mentioned above, as the continuity of the background process
implies continuity in arrival times and guarantees the FIFO-property.
In addition, our setup allows for correlation and is capable of incorporating
both recurrent and non-recurrent congestion. As mentioned above, our second contribution concerns the development of  efficient and close-to-optimal routing algorithms.
The evaluation, through DP,  of the optimal policy being prohibitively slow due to dimensionality issues, we propose the highly efficient
$\edsger$ algorithm. 
Our numerical examples demonstrate that our approach is well capable of  incorporating both recurrent and non-recurrent congestion. Notably, it clearly outperforms deterministic Dijkstra-type algorithms in which the per-arc travel times are replaced by their expected values.
Comparison of the $\edsger$ algorithm to DP shows that $\edsger$ 
is orders of magnitude faster. Indeed, the computational costs of DP grow exponentially with
the network size, while $\edsger$ provides essentially real-time response, 
performing just slightly below optimal.

The background process and SMDP framework are presented in
Section~2. Section~3 considers routing
in this framework, 
describing the optimal routing algorithm, 
heuristic algorithms and speed-up techniques.
Numerical examples to
show the usefulness of the model and to assess the efficiency of the routing
algorithms are presented in Section~4. Section~5 contains concluding remarks.


\section{Markovian Velocity Model}

The first part of this section provides a motivational example for
the proposed framework and the various routing algorithms developed for it.
For illustration we consider a relatively small network, and describe the structure
of a typical background process corresponding to this network.
We briefly assess the expected travel time and run-time of the routing algorithms, 
so as to provide an indication of the gain that can be achieved by our proposed routing scheme.
The second part of the section gives a detailed description of the mathematical
framework for the Markov model that describes the attainable velocities of the vehicles,
and points out how {in the context of this model}  routing
can be interpreted as an SMDP.

\subsection{Motivational example}\label{S21}
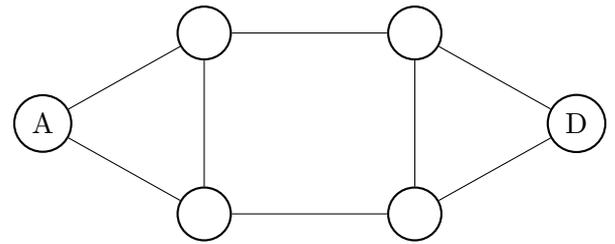
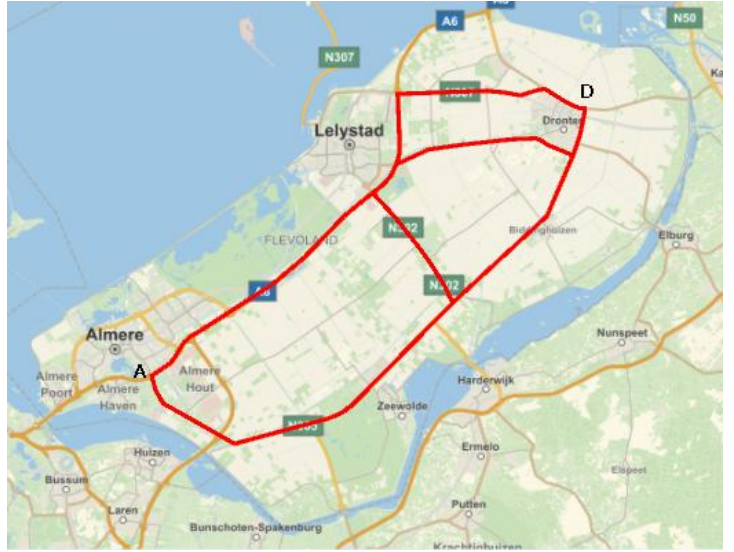
\begin{wrapfigure}{r}{.5\textwidth}
    \begin{minipage}{\linewidth}
    \centering\captionsetup[subfigure]{justification=centering}
    \begin{tikzpicture}
        [
roundnode/.style={circle, draw=black, fill=red!20, thick, minimum size=7mm},
visnode/.style={circle, draw=black, thick, minimum size=7mm},
curnode/.style={circle, draw=red!60, fill=red!20, thick, minimum size=7mm},
squarednode/.style={rectangle, draw=red!60, fill=red!5, very thick, minimum size=5mm},
]
\node[draw=none]    (7)                                                                         {};
\node[draw=none]    (8)     [right=2.5 of 7]                                                    {};
\node[visnode]    (2)     [above=0.7 of 7]    {};
\node[visnode]    (1)     [left= 1.6 of 7]         {A};
\node[visnode]    (4)     [above=0.7 of 8]    {};
\node[visnode]    (3)     [below=0.7 of 7]    {};
\node[visnode]    (5)     [below=0.7 of 8]    {};
\node[visnode]    (6)     [right=1.6 of 8]    {D};
\node[visnode, draw=none]    (7)     [below=0.1 of 8]    {};

\draw[-] (2.south) -- (3.north) node [pos=0.5,left] {};
\draw[-] (4.south) -- (5.north) node [pos=0.5,right] {};
\draw[-] (2.east) -- (4.west) node [pos=0.5,above] {};
\draw[-] (3.east) -- (5.west) node [pos=0.5,below] {};
\draw[-] (1) -- (3) node [pos=0.4,below] {};
\draw[-] (1) -- (2) node [pos=0.4,above] {};
\draw[-] (4) -- (6) node [pos=0.5,above] {};
\draw[-] (5) -- (6) node [pos=0.5,below] {};
    \end{tikzpicture}
    \subcaption{Graphical situation}
    \label{fig:ADa}\par\vfill
    \includegraphics[width=\linewidth]{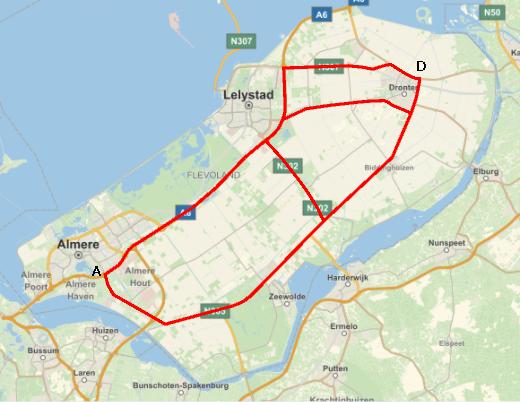}
    \subcaption{Real situation}
    \label{fig:ADb}
\end{minipage}
\caption{Road network Almere-Dronten}\label{fig:AD}
\end{wrapfigure}
Consider the road network of Figure~\ref{fig:AD}{\sc b}, depicting the roads that can be used to travel
from Almere (A) to Dronten (D),
two cities in the Netherlands. As in the rest of this paper, the goal is to minimize the expected travel time.
Figure~\ref{fig:AD}{\sc a} shows the corresponding routing graph in
which each node represents an intersection in Figure~\ref{fig:AD}{\sc b},
and each link represents the road between two of these intersections.
The Dutch government has imposed maximum velocities on the 16 roads
(8 bidirectional arcs, that is) in the network of Figure~\ref{fig:AD}. 
However, due to e.g.\ bad weather conditions and traffic incidents vehicles may
not always be able to drive at these maximum velocities.

We model the variability in velocities by 
a Markovian background process. This background process records the events
that affect these speeds.
In this example, the only events discussed are incidents and rain showers, 
but it can be extended to include various sorts of other events in an obvious manner. 
To model traffic accidents we let
$\{X_{i}(t), t \geq 0\}$ be a Markov process that denotes whether
arc $i$ is blocked by an incident at time $t$. 
Specifically, we choose to set $X_i(t) = 1$
if the arc is not blocked, and $2$ otherwise.
Furthermore, we let $\{Y(t), t \geq 0\}$ be
a Markov process that describes whether it rains at time $t$ or not:
$Y(t) = 1$ indicates it is dry whereas $Y(t) = 2$ indicates it rains.
The Markovian background process can then be written as the vector $B(t) := (Y(t),X_1(t),..,X_{16}(t))$, having  a state space of dimension $2^{17}$.
The state of $B(t)$ describes the velocities
on each of the arcs: we let $v_{i}(s)$ be the velocity on arc $i$ 
when $B(t) = s$.

In the sequel we impose the natural assumption that the processes 
$X_1(t),..,X_{16}(t)$ are independent, which effectively means that the occurrence
of an incident on a specific link has no impact on the occurrence of an incident
on other links. This assumption does not imply that there is no correlation between
travel times on the arcs, as the velocity on an arc depends on the state of $B(t)$.

We let the transition rate matrix for $X_i(t)$, for $i=1,\ldots,16$,
be written as
\begin{align*}
Q_i =
    \begin{bmatrix}
    -\alpha_i & \alpha_i \\
    \beta_i & -\beta_i
    \end{bmatrix},
\end{align*}
with $\alpha_i$ the incident rate and $\beta_i$ the clearance rate
of an incident at arc $i$. 
The $2^{16}\times 2^{16}$-dimensional transition rate matrix 
of the vector process $(X_1(t),\dots,X_{16}(t))$
is now given by $Q_0 := Q_{16} \oplus Q_{15} \oplus \dots \oplus Q_1$; here the operation `$\oplus$' denotes the Kronecker sum \cite{Pease1965},
which, given `$\otimes$' denotes the direct product
and $I_{\dim(C)}$ denotes the identity matrix with the same dimensions
as the matrix $C$, is defined as 
\begin{align}
A \oplus B = A \otimes I_{\dim(B)} 
+ I_{\dim(A)} \otimes B, \label{def:kron}
\end{align}
for two square matrices $A$ and $B$. 
For example, the Kronecker sum for the
two matrices $Q_2$ and $Q_1$ is given by
\begin{align*}
Q_2 \oplus Q_1
= 
    \begin{bmatrix}
    -\alpha_1 - \alpha_2 & \alpha_1 & \alpha_2 & 0 \\
    \beta_1 & -\beta_1 - \alpha_2 & 0 & \alpha_2 \\
    \beta_2 & 0 & -\alpha_1 - \beta_2 & \alpha_1 \\
    0 & \beta_2 & \beta_1 & -\beta_1 - \beta_2
    \end{bmatrix},
\end{align*}
which can be interpreted as the transition rate matrix of $(X_1(t),X_2(t))$ that lives on the state space
$(1,1)$, $(2,1)$, $(1,2)$, $(2,2)$.
Regarding $Y(t)$, we denote the transition rate from $1$ to $2$ as $\lambda$
and from $2$ to $1$ as $\mu$. This means that 
$1/\lambda$ represents the mean time between showers, and $1/\mu$ the mean shower duration. 
Combining the above, the $Q$-matrix for the full vector $B(t)$ can be written
as
\begin{align*}
    Q = \begin{bmatrix}
    Q_0 -\lambda I & \lambda I \\
    \mu I & Q_0 - \mu I
    \end{bmatrix}.
\end{align*}

\begin{remark}{\em {For this example we have assumed that
$X_i(t)$ is independent of $Y(t)$, which informally means that
rain does not affect the
occurrence of incidents. This may sound unnatural, but, importantly, 
our model can be adapted in a straightforward manner to make sure 
that weather conditions do affect the rate of incidents. 
The general setup will be introduced in full detail in the next subsection.}}$\hfill\Diamond$
\end{remark}

In the model introduced above the travel time on 
any arc is completely determined by the proposed velocity dynamics.
This means that for any given arc, if the state of the background process when leaving the origin node is given, 
the expected travel time to the destination node 
can be computed {(jointly with the state of the background process when arriving at the destination node).} 
This in principle allows the computation of the optimal routing policy.
{Below we will argue that this optimal policy can be
found relying on dynamic programming (DP) methods}. 
{DP-based methods}, however, are typically prohibitively
computationally demanding
in case  the underlying state space is large (the {\it curse of dimensionality}). Focusing on our specific routing objective, it is clear that such computational issues will arise, since,
as we saw above, the number of background states is large,
even in a small network.

The $\edsger$ algorithm, that will be presented in Section~\ref{DRA}, succeeds in overcoming the high computation complexity
of DP-based methods.
The idea behind $\edsger$ 
is inspired by
the A$^\star$-algorithm, 
a speed-up version of Dijkstra's algorithm. 
We assess $\edsger$ on two criteria, namely:
\begin{itemize}
    \item[$\circ$] {\it distance-to-optimality}, i.e., the
    difference between the value of the objective function and the optimal value.
    As the results of $\edsger$ are just a marginal amount below the optimal value,
    $\edsger$ is close-to-optimal.
    \item[$\circ$] {\it efficiency}, i.e., the run-time of the algorithm.
    $\edsger$ is highly efficient as it is sufficiently fast to allow it being used in real-time.
\end{itemize}

To provide an impression of the achievable performance of $\edsger$, we conducted an experiment on the network presented in Figure~\ref{fig:AD}. We define  the $17$-dimensional
background process $B(t)$
as  above. We chose
$\alpha_i = 0.1 \:{\rm h}^{-1}, \beta_i = 2 \:{\rm h}^{-1}$ for $i = 1,\dots,16$ and $\lambda = \mu = 0.25 \:{\rm h}^{-1}$, meaning that
the expected time between two incidents is $10$ hours and the expected clearance
time of an incident $30$ minutes, whereas both 
the expected duration of a rain shower and the expected duration of a dry period equal $4$ hours.
{These parameter values were chosen merely for illustrative purposes; in all later experiments we base ourselves on historical data.}
In  case that it does not rain and there is no incident on the arc  we let the vehicle speed be 120 km/h if  there is no incident on the directly adjacent arcs, and 100 km/h otherwise. Regarding the case that it does not rain and
there is an incident on the arc, we let the vehicle speed
be 50 km/h if there is no incident on the directly adjacent
arcs, and 20 km/h otherwise. In case it does rain the velocities on this arc,
in the four situations discussed above, are 100 km/h, 80 km/h, 20 km/h and 10 km/h, respectively.

\begin{table}[h]
\begin{tabular}{|c||c|c||c|}
\hline
    & Average Expected Travel Time & Weighted Average Expected Travel Time & Run-time (sec)\\
   DS & $0.510$ & $0.433$ & $1.35 \cdot 10^{-3}$ \\
    $\edsger$  & $0.443$ & $0.410$ & $1.13 \cdot 10^{-2}$ \\
    {DP} & $0.442$ & $0.410$ & $4.26 \cdot 10^{-1}$\\
    \hline
\end{tabular}
    \caption{Results example Almere-Dronten}
    \label{tab:AD}
\end{table}

The results pertaining to this example are shown in Table ~\ref{tab:AD}.
The performance of the $\edsger$ algorithm is compared
with a deterministic Dijkstra-type algorithm that assumes that a vehicle can always
drive at the maximum speed level and thus does not take any stochasticity
into account (`DS', being an abbreviation of `Deterministic Static').
We in addition implemented a competitive {DP} algorithm for determining the
policy that minimizes the expected travel time 
{(`DP')}.
\begin{itemize}
     \item[$\circ$]
    
The first column  shows the expected travel time
from A to D under the different policies, averaged (evenly) over
all possible initial background states.
 \item[$\circ$] The second column
provides a weighted average of the expectations corresponding to 
all possible initial background states; for a given initial background state,
the weight is chosen 
equal to its limiting probability.
 \item[$\circ$]
The last column contains the run-time of the three algorithms.

\end{itemize}
A first conclusion is that
the objective function achieved by $\edsger$ is close to its minimal value (as provided by DP). 
{Secondly}, even in this small network, the run-time
of $\edsger$ is significantly lower than the run-time
of the competitive 
{DP}
algorithm.
Numerical experiments later in this paper will show
that the run-time of the competitive 
{DP}
algorithm grows exponentially
with the network size, whereas the run-time of the $\edsger$ 
algorithm remains essentially real-time.
We in addition conclude that ignoring the stochasticity, as is done by DS, 
leads to a fast but 
far from optimal
algorithm; note in 
particular that the corresponding objective function is substantially 
higher than its minimal value (as provided by DP).

\subsection{Mathematical Framework}
After the motivating example of the previous subsection, 
we now formally introduce the general mathematical framework that we will be working with.
As before, we consider the objective of minimizing the expected
travel time between a given origin-destination (OD)
pair in the road network. 
Let $G = (N,A)$ be a graph representing the road network, 
with $N$ and $A$ the set of nodes in $G$, and 
the set of directed arcs in $G$, respectively.
The set $N$ represents the intersections in the road
network, whereas the set $A$ represents the roads
that connect these intersections, implying $k\ell \in A$
only if there is a (direct) road in the network between the intersections that are 
labeled by $k$ and $\ell$.
Let $d_{k\ell}$ be the length of arc $k\ell$.

A realistic feature of our setup is that when this arc is traveled by a vehicle, the
speed of this vehicle is not necessarily constant. Indeed, it
may vary between finitely many values, related to e.g.\ the occurrence of incidents and weather conditions.
To deal with these changing velocities, we introduce a
Markovian background process on the arcs $A$, of which the background process discussed in Section~\ref{S21} is a special case.
As we will argue below, the use of this stochastic process will allow us
to incorporate random effects (corresponding to non-recurrent congestion, that is) as well as (near-)deterministic patterns
in the attainable speeds on the arcs (corresponding to recurrent congestion).
We assume that the state of the background process
and the corresponding traffic velocities are available while traveling. Indeed, travelers are able to adapt
the chosen route while driving,
based on the information available. 
In the sequel we will phrase this dynamic stochastic shortest
path problem with information on velocity levels
in terms of a finite semi-Markov decision model.

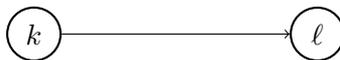
\begin{figure}[h]
\begin{center}
\begin{tikzpicture}[
roundnode/.style={circle, draw=black, thick, minimum size=7mm},
visnode/.style={circle, draw=black, thick, minimum size=7mm},
curnode/.style={circle, draw=red!60, fill=red!20, thick, minimum size=7mm},
squarednode/.style={rectangle, draw=red!60, fill=red!5, very thick, minimum size=5mm},
]
\node[roundnode]    (1)     {$k$};
\node[roundnode]    (2)     [right= 3 of 1]         {$\ell$};

\draw[->] (1.east) -- (2.west) node [pos=0.5,above] {};
\end{tikzpicture}
\caption{Single arc network}
\label{fig:singlearc}
\end{center}
\end{figure}

{{\sc Background process}} \\
Before presenting the full Markov-modulated environmental process,
for expositional reasons we will first show examples of background processes on small networks 
that fit into our framework and gradually
extend this setup. 
First consider the single-arc network in Figure~\ref{fig:singlearc}.
Define $\{X_{k\ell}(t), t \geqslant 0\}$ as the continuous-time 
Markovian background process
on this arc such that
\begin{align*}
    X_{k\ell}(t) &= 
    \begin{cases}
    1 &\quad \quad \text{no incident on arc $k\ell$ at time $t$} \\
    2 &\quad \quad \text{otherwise},
    \end{cases}
\end{align*}
similar to the modeling of incidents in the motivational example of
the preceding subsection.
The process $X_{k\ell}(t)$ provides information on possible
events on $k\ell$ and determines 
the velocity of a vehicle traveling on $k\ell$: the speed 
at time $t$ equals $v_{k\ell}(i) $
if $X_{k\ell}(t) = i$, for $i=1, 2$. 
For technical reasons it is throughout assumed that all these velocities are strictly positive, 
but this is, in practical terms, not a restriction as they are allowed to have arbitrarily small values.
The speed at which vehicles can travel on this arc is now completely described
through the background process $X_{k\ell}(t)$ and the corresponding $v_{k\ell}(i)$.

In the above example there are two possible speeds, but note that
in practice one could work with more than two
levels; one could e.g.\ think of the period between the clearance
of an incident and the time the free-flow speed can be attained again.
These dynamics can be incorporated by allowing the state space of the process 
$X_{k\ell}(t)$ to consist of more than two states.
In general we let this state space be denoted by 
$\mathcal{S}_{k\ell} = \{1,\dots,n_{k\ell}\}$. When $X_{k\ell}(t) = s \in \mathcal{S}_{k\ell}$,
the speed at which vehicles are moving on arc $k\ell$
is $v_{k\ell}(s)$.
We 
use the notation $Q_{k\ell}$ to denote the corresponding transition rate matrix
of dimension $n_{k\ell} \times n_{k\ell}$.

\begin{figure}[h]
\begin{center}
\begin{tikzpicture}[
roundnode/.style={circle, draw=black, thick, minimum size=7mm},
visnode/.style={circle, draw=black, thick, minimum size=7mm},
curnode/.style={circle, draw=red!60, fill=red!20, thick, minimum size=7mm},
squarednode/.style={rectangle, draw=red!60, fill=red!5, very thick, minimum size=5mm},
]
\node[roundnode]    (1)     {$k$};
\node[roundnode]    (2)     [right= 2 of 1]         {$\ell$};
\node[roundnode]    (3)     [right= 2 of 2]         {$m$};

\draw[->] (1.east) -- (2.west) node [pos=0.5,above] {};
\draw[->] (2.east) -- (3.west) node [pos=0.5,above] {};
\end{tikzpicture}
\caption{Two-arc network}
\label{fig:twoarc}
\end{center}
\end{figure}
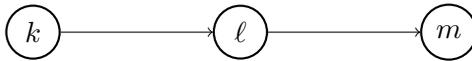

We now extend the network with an arc $\ell m$ (Figure~\ref{fig:twoarc}), and
let $B(t) = (X_{k\ell}(t),X_{\ell m}(t))$ be the Markovian background process
on this two-arc network.
Here $\{X_{k\ell}(t),t\geqslant 0\}$ and $\{X_{\ell m}(t),t\geqslant 0\}$ are 
independent
Markovian background processes on the arcs $k\ell$
and $\ell m$, respectively.
As in the single-arc case we denote the
state space 
of $X_{k\ell}(t)$ by $\mathcal{S}_{k\ell} = \{1,\dots,n_{k\ell}\}$
and the transition rate matrix as $Q_{k\ell}$. 
We can equivalently define $\mathcal{S}_{\ell m}$ and
$Q_{\ell m}$ for the process $X_{\ell m}(t)$.
With $X_{k\ell}(t)$ and $X_{\ell m}(t)$ being independent,
the transition rate matrix $Q$ of $B(t)$
is of the form
\begin{align}
    Q = Q_{k\ell} \oplus Q_{\ell m},
\end{align}
where `$\oplus$' is the Kronecker sum that was defined
in (\ref{def:kron}).
If $B(t)$ is in state
$s \in \mathcal{I} := \mathcal{S}_{k\ell} \times \mathcal{S}_{\ell m}$,
we let the speed at which vehicles are moving 
be $v_{k\ell}(s) $ on arc $k\ell$ and $v_{\ell m}(s) $
on arc $\ell m$.
Importantly, this means that the velocity on each of the two arcs is allowed to depend on both $X_{k\ell}(t)$
and $X_{\ell m}(t)$. 
Hence, this way of modeling 
defines an implicit dependence between the velocities on the two arcs.

The presented setup can easily be extended by
adding more arcs to the network.
Consider the network $G = (N,A)$
and label the directed edges 
in $A$ by $k_1\ell_1,k_2\ell_2,\dots,k_n\ell_n$ with $n:=|A|$.
The Markovian background process 
can now be written as $B(t) = (X_{k_1\ell_1}(t),\dots,X_{k_n\ell_n}(t))$.
Here $\{X_{k_j\ell_j}(t), t \geqslant 0\}$, with $j\in\{1,\ldots,n\}$, is a
continuous-time Markov
process with state space
${\mathcal S}_{k_j\ell_j} = \{1,\dots,n_{k_j\ell_j}\}$ and transition rate matrix
$Q_{k_j\ell_j}$, such that $X_{k_j\ell_j}(t)$ is the congestion level
at arc $k_j\ell_j \in A$ at time $t$. Assuming that
$X_{k_i\ell_i}(t), X_{k_j\ell_j}(t)$ evolve independently for $i \neq j$,
the transition rate matrix $Q$ of $B(t)$
is of the form
\begin{align*}
    Q = Q_{k_n\ell_n} \oplus Q_{k_{n-1}\ell_{n-1}} \oplus \cdots 
    \oplus Q_{k_2\ell_2} \oplus Q_{k_1\ell_1}.
\end{align*}
Dependence between the arcs can again be realized by 
letting the speed on each of the arcs depend
on the state of the full vector $B(t)$:
when $B(t)$ is in state $s \in \mathcal{I} = \mathcal{S}_{k_1\ell_1} \times 
\dots \times \mathcal{S}_{k_n\ell_n}$ the speed
at which vehicles are moving on arc $k_j\ell_j$ is
$v_{k_j\ell_j}(s)$.
The dynamics of the vehicles moving on the network are now completely defined
through the background process $B(t)$ and the corresponding
velocity levels. For instance, if $B(u)=s$, then
the travel time
on edge $k_i\ell_i$ for traveling a distance
$d \in [0,d_{k_i\ell_i}]$ when leaving node $k_i$ at
time $u$ is distributed as
$\tau_{k_i\ell_i}^{s}(d)$, with
\[\tau_{k_i\ell_i}^{s}(d):= \min\left\{t\geqslant 0: \left.
\int_{0}^{t} v_{k_i\ell_i}( B(u))\,{\rm d}u
\geqslant d \:\right|\, B(0) = s\right\},\]
where we
use the fact that the travel distance $d$ 
can be computed by an integral over the travel speed.
For transparency of notation we abbreviate 
\[\tau_{k_i\ell_i}^{s}(d_{k_i\ell_i})\equiv
\tau_{k_i\ell_i}^{s}\] to denote the
travel time for arc $k_i\ell_i$ when leaving at time $u$ from node $k_i$.

There is global correlation in the above setup, due to the fact that the speed that a vehicle on a given arc can drive at  depends on the 
full background process
$B(t)$, containing information on the congestion levels of {\it all} arcs.
Previous research has shown that correlation
due to non-recurrent congestion is primarily local; the congestion
level of an arc mainly affects the attainable 
velocities at upstream arcs that are within a certain distance
of the incident arc \cite{Guo2019, Priambodo2020}. We model this property through the following assumption.
\begin{assumption}[local-$r$-correlation]
Denote with $A_{k\ell}^r$ the set of arcs 
that are at most $r$ arcs away from arc $k\ell$ \emph{(}including
$k\ell$ itself\emph{)}.
We assume the velocity on arc $k\ell$ to depend only 
on the congestion levels of the arcs at most $r$ arcs away:
$v_{k\ell}(s) = v_{k\ell}(s')$ if $s_i = s_i'$ for all
$i \in A_{k\ell}^r$.
\end{assumption}

Before continuing to the
SMDP routing framework induced by our proposed model, 
the next remark discusses a possible extension of our model, which was already 
noted in the motivational example.

\begin{remark}{\em 
{Suppose there is} a Markovian process $Y(t)$ that prompts global correlation
in the network, e.g.\ weather conditions that affect
the velocities on {\it all} arcs in the network.
This can be incorporated into our framework by
expanding the process $B(t)$ with the process $Y(t)$.
We let $B(t) = (Y(t),X_{k_1\ell_1}(t),\dots,X_{k_n\ell_n}(t))$
be such that, given the state of $Y(t)$,
the future evolution of 
$X_{k_i\ell_i}(t)$ and $X_{k_j\ell_j}(t)$ are independent for $i \neq j$.
The transition rate matrix for the case $Y(t) = s$ is then given by
\begin{align}
    Q^s = Q_{k_n\ell_n}^s \oplus Q_{k_{n-1}\ell_{n-1}}^s \oplus 
    \dots \oplus Q_{k_2\ell_2}^s \oplus Q_{k_1\ell_1}^s, \label{kronwithy2}
\end{align}
with $Q_{k_j\ell_j}^s$ the transition rate matrix
for the congestion levels at arc $k_j\ell_j$ given 
$Y(t) = s$. In this way our model allows for global dependence when
information on processes that affect
velocities in the whole network is available.
In the remainder of this paper we will 
mainly focus on the case without global correlation,
but, by using transition matrices of the type (\ref{kronwithy2}), the results of this paper can be extended
to the case of a known global dependence process.}\hfill$\Diamond$
\end{remark}

\begin{remark}
\label{rem:erlang}
{\em 
By the use of phase-type distributions the proposed
Markov model for congestion
is not limited to exponentially distributed 
times between states of congestion \cite{Asmussen2003}.
Importantly, phase-type distributions can model random quantities that are less variable than the exponential distribution (e.g.\ by using Erlang distributions)
as well as random quantities that are more variable than the exponential distribution (e.g.\ by using hyperexponential distributions). In particular, for highly predictable events (think of recurring events, such as rush hours) Erlang distributions with a high number of phases, thus leading to a low variance, are well suited. 
}\hfill$\Diamond$
\end{remark}

{{\sc Semi-Markov Decision Process}}\\
The objective of this subsection is to phrase our optimization problem in terms of an SMDP.
Above we introduced our stochastic process modeling the  velocity dynamics.
Our setup
is somewhat in the spirit of the one used in 
the work of Psaraftis and Tsitsiklis \cite{Psaraftis1993}, 
Kim, Lewis and White
\cite{Kim2005a, Kim2005b} and Sever et al.\ \cite{Sever2013}, 
who use an MDP
to capture the stochasticity in travel times.
Importantly, as opposed to such earlier MDP-based approaches, 
we impose a {\it  continuous-time}, rather than discrete-time, stochastic process on the arc speeds. 
{Although at the expense of additional computational
complexity, this 
has several advantages}:
\begin{itemize}
    \item [$\circ$]When a vehicle travels an arc, 
    changing conditions on this arc are taken into account.
    In the conventional MDP framework it is assumed that if a vehicle
    is at an intersection, travel times on the attached arcs are known
    and choices are based on these known travel times. 
    Conditions on these attached
    arcs may however change while driving, 
    affecting the arrival time.
    \item [$\circ$]Travel times are a continuous function of the departure epoch.
    Therefore the consistency- or FIFO-property \cite{Wellman1995},
    is naturally satisfied (see Appendix A for the proof):
    \begin{proposition}
    \label{prop1}
    Let $k\ell \in A$ and denote with $\tau_{k\ell}^{t}(d)$ the travel time on arc
    $k\ell$ for traveling a distance $d \in [0,d_{k\ell}]$ when starting
    at $t \geqslant 0$. Then 
    $t_1 \leqslant t_2$ implies that $t_1 + \tau_{k\ell}^{t_1}(d)
    \leqslant t_2 + \tau_{k\ell}^{t_2}(d)$.
    \end{proposition}
    \item [$\circ$]In both our framework and the MDP setting 
    the Markov processes on the links are assumed independent.
    In the MDP framework this implies independence 
    between the travel times on the links. 
    By imposing a background process
    affecting the arc velocities we however
    provide the opportunity to incorporate
    dependence between the travel times on the links.
\end{itemize}


Routing in our framework can be analyzed using 
a finite semi-Markov decision process (SMDP).
The state of the SMDP consists of the location of the
vehicle combined with the state of the background process.
In a practical context, the first could be tracked
by GPS, whereas the latter can be provided by an ITS.
Decision epochs are the arrival times at the intersections,
in the sense that at these epochs one has the opportunity to adapt the route.
Denote at decision-epoch $i$ the state of the SMDP
as $(K_i, B(t_i))$, with $K_i \in N$
the label of the current intersection
and $B(t_i)$ the state of the background process
at the current time $t_i$. 
The following (reasonable) assumptions are made:
\begin{itemize}
    \item [$\circ$]Waiting at a node is not allowed. 
    Note that this assumption does not
    affect an optimal policy, as the FIFO property implies that waiting
    is never advantageous.
    \item [$\circ$]Information on the congestion levels on all
    arcs is available at all times.
    \item [$\circ$]There is at least one path
    connecting every node with the destination node.
    If this assumption is not fulfilled, we have
    a non-connected graph for which we can only construct routing policies
    for OD-pairs within connected components.
\end{itemize}

At decision epoch $i$, 
the goal is to choose a neighbor $k_{i+1}$ of $k_i$
such that the expected travel time from $k_i$ 
to the destination $k^\star$ is minimized.
A policy matrix assigns a successor node for each node and each current state 
of the background process (i.e., each pair $(k,s)$
with $k \in N, s \in \mathcal{I}$).
Given destination $k^\star$, the expected travel time under a 
policy $\pi$ and initial state $(K_0,B(t_0))$ 
is given by
\begin{align}
    J^{\pi}(K_0,B(t_0)) = \mathbb{E}\left[\sum_{i = 0}^{I^\star-1}
    \tau_{K_iK_{i+1}}^{B(t_i)} \,\bigg|\, K_{i+1}=\pi(K_i, B(t_i)), \, \enspace 
    {i = 1,\dots,I^\star\!-\!1}\right] \label{eq:Vvol},
\end{align}
in which $I^\star$ denotes the number of decision epochs 
until the destination $k^\star$ is reached, $t_i$ the
time of decision epoch $i$,
and $\pi(K_i,B(t_i))$ the action under policy
$\pi$ given the pair $(K_i,B(t_i))$.
The optimal policy is now a policy $\pi^\star$ such that 
\begin{align*}
    \pi^\star = \argmin_{\pi \in \Pi} J^{\pi}(K_0,B(t_0)),
\end{align*}
with $\Pi$ denoting the set of admissible policies (i.e., 
policies for which the probability of reaching
the destination is one for all initial states).
For every non-admissible policy $\pi$ set $J^\pi(k,s) := \infty$
for all $(k,s) \in N \times \mathcal{I}$.

We conclude this section by deriving expressions for both the expected travel time
and the transition probabilities pertaining to a single arc. The resulting quantities are building
blocks of the routing policies that will be discussed in the next section.
To this end, define, for $d \in [0,d_{k\ell}]$
and $s,s' \in \mathcal{I}$, the expected travel time on $k\ell$:
\begin{align*}
\phi_s(d \,|\, k,\ell) &:= \mathbb{E}[\tau_{k\ell}^s(d)] ,
\hspace{10mm} \Phi(d \,|\, k,\ell) := (\phi_s(d\,|\,k,\ell))_{s \in \mathcal{I}}
\end{align*}
In addition, we will work with the Laplace-Stieltjes transform of the travel time on edge $k\ell$
intersected with the event that upon completion the 
background state is $s'$:
\begin{align*}
    \psi_{ss'}(d \,|\, \alpha, k,\ell) &:= \mathbb{E}[e^{-\alpha\tau_{k\ell}^s(d)}1\{
B(\tau^s_{k\ell}(d))=s'\}],
\hspace{10mm} \Psi(d \,|\, \alpha, k, \ell) := (\psi_{ss'}(d\,|\,\alpha,k,\ell))_{(s,s') \in \mathcal{I} \times \mathcal{I}}.
\end{align*}
The objects introduced above can be evaluated by conditioning
on a possible jump of the background process in a small time-interval.
As shown in the following theorem, it leads to a system of linear differential equations,
whose solution can be given in terms of matrix
exponentials. A proof is given
in Appendix B.

\begin{theorem}
\label{Theorem1}
Given a graph $G = (N,A)$ with a pair of nodes $k,\ell \in N$ and a distance $d\geqslant 0$,
it holds that
\begin{align*}
    V_{k\ell}\, \Phi'(d \,|\, k,\ell) &= {\boldsymbol 1} + Q\,\Phi(d \,|\, k,\ell), \\
    V_{k\ell}\, \Psi'(d \,|\, \alpha,k,\ell) &= (Q - \alpha I)\,\Psi(d \,|\, \alpha,k,\ell),
\end{align*}
with $V_{k\ell} := {\rm diag}\,\{(v_{k\ell}(s))_{s \in \mathcal{I}}\}$
and ${\boldsymbol 1}$ a $|\mathcal{I}|$-dimensional column vector of ones.
A solution for this system of linear differential equations can be written as
\begin{align}
    \Phi(d \,|\, k,\ell)
    &= 
    \exp\left\{
    \begin{bmatrix}[1.3]
    dV_{k\ell}^{-1}Q & dV_{k\ell}^{-1}{\boldsymbol 1} \\ {\boldsymbol 0}^\top & 0 \end{bmatrix}\right\}
    \begin{bmatrix}
    {\boldsymbol 0} \\ 1
    \end{bmatrix} \label{eq:solphi} \\
    \Psi(d \,|\, \alpha,k,\ell) &=
    \exp\{d\,V_{k\ell}^{-1}(Q - \alpha I)\} \label{eq:solpsi}
\end{align}
with ${\boldsymbol 0}$ an $|\mathcal{I}|$-dimensional
column vector of zeroes.
\end{theorem}
Due to our requirement that all the $v_{k\ell}(s)$ are positive,
the matrix $V_{k\ell}^{-1}$ is well-defined. 
With the above result the expected travel time on an arc
can be directly computed using the expression for
$\Phi(d \,|\, k,\ell)$. 
An expression for 
the transition probabilities can be found by 
{setting} $\alpha$ 
{equal to} 0 in (\ref{eq:solpsi}):
\begin{corollary}
\label{cor:transprob}
For a pair of nodes $k,\ell\in N$ and pair of background states $s,s'$,
\begin{align}
    \mathbb{P}(B(\tau^s_{k\ell}) = s') = \big[\exp\big(d_{k\ell} V_{k\ell}^{-1} Q \big)\big]_{s,s'} \label{eq:trans}
\end{align}
\end{corollary}
\begin{remark}
\label{rem:suff}
{\em 
The upper left $|\mathcal{I}| \times |\mathcal{I}|$-block of
the matrix exponential in (\ref{eq:solphi})
is the matrix exponential of $dV_{k\ell}^{-1}Q$. 
Thus to compute both the expectation
$\mathbb{E}[\tau^s_{k\ell}]$
and the transition probabilities $\mathbb{P}(B(\tau^s_{k\ell}) = s')$, 
computation of the matrix exponential given in {(\ref{eq:solphi})}
suffices.}\hfill$\Diamond$
\end{remark}

\section{Dynamic Routing Algorithms}\label{DRA}

Recall that our objective is to minimize the expected travel time
for a given OD-pair in a road network in which 
a vehicle may experience changing velocities and in which one knows the current state of the driving background process $B(t)$. After having provided background 
on the optimal policy resulting from DP  (Section~3.1), 
we  present 
two Dijkstra-like shortest path algorithms
that aim to output a (near-)optimal routing policy.
{These algorithms are dynamic shortest path algorithms, in that at every new
intersection a shortest path algorithm is run again}.
The analysis of Section~3.2 reveals that the first of these two algorithms, which we have called {\sc edsger},  suffers from the curse
of dimensionality. 
The second presented algorithm, $\edsger$, as 
{introduced} in Section~3.3, overcomes
this drawback, and offers real-time response.
Useful implementation details, for these two
algorithms as well as for the DP-based approach, will be provided in Section~3.4.
We conclude with Section~3.5, which presents speed-up techniques that reduce
the state space and/or network size, to be used to further decrease the computational costs.


\subsection{Optimal Policy by Dynamic Programming}
In this subsection we present an algorithm that outputs an optimal routing policy. As we will see in our numerical examples, this method is prohibitively slow; in this paper we mainly use it as a benchmark for our routing algorithms {\sc edsger} and $\edsger$. 
{We will
argue below that} optimal policies can be characterized
as solutions of a set of Bellman optimality equations
\cite{Bellman1957a, Bellman1957b, bertsekas1976, Bertsekas2005, Puterman1994, Ross1983, Sutton2018}.
DP methods, competitive algorithms that solve these optimality equations,
have high computational costs, so that they are not particularly suitable for real-time routing purposes.
Details on the implementation of the competitive DP algorithm are 
provided in Section~3.4. 

We proceed by arguing 
that algorithms outputting an optimal routing policy suffer
from the curse of dimensionality.
In case of a MDP,
G\"uner et al.\ \cite{Guner2012}, Kim et al.\ \cite{Kim2005a, Kim2005b}
and Sever et al.\ \cite{Sever2013} have characterized an optimal policy
to satisfy the Bellman optimality equations.
Importantly, this property carries over
to our framework, as can be seen as follows.
{First note that we have a destination node $k^\star$
which is cost-free and absorbing, i.e., 
$k^\star = \pi(k^\star, s)$ and 
$J^\pi(k^\star, s) = 0$ for all $s \in \mathcal{I}.$
Furthermore,} denoting $k_1 := \pi(k_0,B(t_0))$ as the next node
under an admissible policy $\pi$ given the state $(k_0, B(t_0))$ with $k_0 \in N 
\setminus k^\star$ and $t_0 \in \mathbb{R}$, (\ref{eq:Vvol}) implies
\begin{align*}
    J^{\pi}(k_0,B(t_0)) &= \mathbb{E}\big[\tau_{k_0k_1}^{B(t_0)}\big]
    + \sum_{s' \in \mathcal{I}} \mathbb{P}\big(B(\tau_{k_0k_1}^{B(t_0)}) = s'\big)
    \,J^{\pi}(k_1,s').
\end{align*}
Thus, any admissible policy satisfies a set of Bellman equations.
Denote with $\NB(k)$ the set of neighbors of a node $k \in N$, i.e., 
$\ell \in \NB(k)$ only if $k\ell \in A$ (note that
$k \notin \NB(k)$ as {it is not allowed to wait
at a node}). 
An optimal policy $\pi^\star $ can then be
given through the Bellman optimality equations:
\begin{align}
    \pi^\star(k_0,B(t_0)) = \argmin_{k_1 \in \NB(k_0)}\left\{\mathbb{E}\big[\tau_{k_0k_1}^{B(t_0)}\big]
    + \sum_{s' \in \mathcal{I}} \mathbb{P}\big(B(\tau_{k_0k_1}^{B(t_0)}) = s'\big)\,
    J^{\pi^\star }(k_1,s')\right\}. \label{eq:Vuiteen}
\end{align}
The expected travel times and transition 
probabilities in (\ref{eq:Vuiteen}) can be computed
relying on the expressions derived in the previous section (in particular Equations  (\ref{eq:solphi}) and (\ref{eq:trans})).

Algorithms that solve the Bellman optimality equations have 
received substantial attention in the literature.
Examples of frequently used solution methods
are 
{
the dynamic programming (DP)
methods (e.g. value iteration and policy iteration)
and linear programming}.
It should be realized that the {size} of our state space may explode: 
even in a simple setup in which every link has only two possible background states, a
network with $|A|$ links leads to a state space of size $2^{|A|}$.
We will therefore employ
value iteration (VI) to solve the optimality
equations \cite{Sutton2018}.
For large state spaces the VI algorithm
does however suffer from high complexity:
\begin{itemize}
    \item [$\circ$] 
        in every iteration of the algorithm the 
        right hand side of (\ref{eq:Vuiteen}) is computed for every 
        possible $B(t_0) \in \mathcal{I}$, so that the complexity of the value
        iteration algorithm contains a 
        term $|\mathcal{I}|^2$;
    \item [$\circ$]
        computation of the expectations and the transition probabilities, which requires the evaluation of a matrix exponential as in (\ref{eq:trans}), will be intractable when the transition rate matrix $Q$ is large;
    \item [$\circ$]a large policy matrix will lead to
        high {memory} costs.
\end{itemize}
{DP suffers from the curse of dimensionality; expanding
the network leads to an exponential increase in the size of the state space 
and corresponding exponential increase in computational costs.}

Based on the above, it can be concluded that VI  
has limited potential in  practical settings.
We will therefore introduce two dynamically-used
Dijkstra-like shortest path algorithms, {\sc edsger} and $\edsger$; 
here `dynamically-used' means that the shortest path 
algorithm is run at any decision epoch along the way.
{\sc edsger} provides near-optimal solutions,
but, similar to VI, suffers from the curse of dimensionality.
$\edsger$ is an adaptation of {\sc edsger}
that uses the local-$r$-correlation assumption to overcome
the curse of dimensionality.
Importantly, the algorithm offers real-time response,
while still providing near-optimal solutions.

\subsection{\edsgernonstar algorithm}
The {\sc edsger} algorithm (where {\sc edsger} is an 
abbreviation of `Expected Delay in Stochastic Graphs
with Efficient Routing') can be seen as a dynamically-used stochastic
version of the A$^\star$-algorithm.
More concretely, the method works as follows:
\begin{itemize}
    \item[$\circ$] At every decision epoch (every
arrival at a node, that is)
we {use a stochastic shortest path algorithm to} 
identify a path with lowest expected 
cost, from the current
node to the destination, provided that one is {\it not} allowed 
to change the route along the way.
The resulting policy $\pi_{\text{E}}$ is to travel to the
next node along this path.
Thus, given that a vehicle is at node $k_0 \in N$ at time $t_0$,
{\sc edsger} uses 
{a shortest path}
algorithm with input ($k_0, B(t_0)$)
to output a path from $k_0$ to $k^\star$.
Then, denoting $k_1$ for the first node in this path,
we have that $\pi_{\text{E}}(k_0, B(t_0)) = k_1$.
\item[$\circ$]
This procedure is then `dynamically-used' in the sense that it is  
repeated when arriving at $k_1$, so as to identify
$k_2 := \pi_{\text{E}}(k_1,B(t_1))$. Given the fact that the state of the background process may have changed while traveling from $k_0$ to $k_1$, this next node may differ from the one that was selected at $k_0$. The procedure thus exploits the information currently available. 
\item[$\circ$] One proceeds along these lines until the destination node $k^\star$ has been reached.
\end{itemize}

A first important remark is that this procedure is not necessarily minimizing the expected delay, i.e., 
the resulting route does not always coincide with the one generated by the DP-based algorithm discussed earlier. 
We note that one reason for this potential loss
is the fact that {\sc edsger} uses a stochastic A$^\star$-like shortest path
algorithm. In a network with {\it stochastic} travel times
the notion that a subpath of a shortest path is a
shortest path, a property on which Dijkstra and
A$^\star$ are built, {does not generally hold}, {so that the output is not guaranteed to be the
shortest path in expectation.}
An example (in 
our context) of this finding, which is
often 
attributed to Hall \cite{Hall1986}, is given in the following example.

\begin{figure}[h]
\begin{center}
\begin{tikzpicture}[
roundnode/.style={circle, draw=black, thick, minimum size=7mm},
extranode/.style={circle, draw=black, thick, minimum size=4mm},
nonvisnode1/.style={circle, draw=white, thick, minimum size=1mm},
nonvisnode2/.style={rectangle, draw=white, thick, minimum width=17mm, minimum height = 10mm},
visnode/.style={circle, draw=black, thick, minimum size=7mm},
curnode/.style={circle, draw=red!60, fill=red!20, thick, minimum size=7mm},
squarednode/.style={rectangle, draw=red!60, fill=red!5, very thick, minimum size=5mm},
]
\node[roundnode]    (1)     {$k_0$};
\node[roundnode]    (2)     [right= 2 of 1]         {$k_1$};
\node[roundnode]    (3)     [right= 3 of 2]         {$k^\star$};
\node[nonvisnode2]  (5)     [right= 2 of 3]         {$(k_0k_1)_2$:};
\node[nonvisnode2]  (6)     [below= 0.1 of 5]       {$k_1k^\star$:};
\node[extranode]    (7)     [right= 1 of 5]         {1};
\node[extranode]    (8)     [right= 1 of 7]         {2};
\node[extranode]    (9)     [right= 1 of 6]         {1};
\node[extranode]    (10)    [right= 1 of 9]         {2};
\node[extranode]    (11)    [right= 1 of 10]        {3};
\node[nonvisnode2]  (12)    [above= 0.1 of 5]       {$(k_0k_1)_1$:};
\node[extranode]    (13)    [right= 1 of 12]        {1};

\draw[->, bend angle = 45, bend left] (1) to node [above,pos=0.5] {$v\big((1,\cdot,\cdot)\big) = 60$} (2);
\draw[->, bend angle = 45, bend right] (1) to node [text width = 3cm, pos=0.5,below] {$v\big((\cdot,1,\cdot)\big) = 100$ $v\big((\cdot,2,\cdot)\big) = 10$} (2);
\draw[->] (2.east) -- (3.west) node [text width = 3cm, pos=0.5,above] {$v\big((\cdot,\cdot,1)\big) = 100$ $v\big((\cdot,\cdot,2)\big) = 80$
$v\big((\cdot,\cdot,3)\big) = 2$};
\draw[->, bend angle = 35, bend left] (7) to node [above,pos=0.5] {$1$} (8);
\draw[->, bend angle = 35, bend left] (8) to node [below,pos=0.5] {$1$} (7);
\draw[->] (9) to node [below,pos=0.5] {$1$} (10);
\draw[->] (10) to node [below,pos=0.5] {$1$} (11);
\end{tikzpicture}
\caption{Example in which {\sc edsger} is not optimal, as optimal subpath is not in optimal path.}
\label{fig:subpath}
\end{center}
\end{figure}
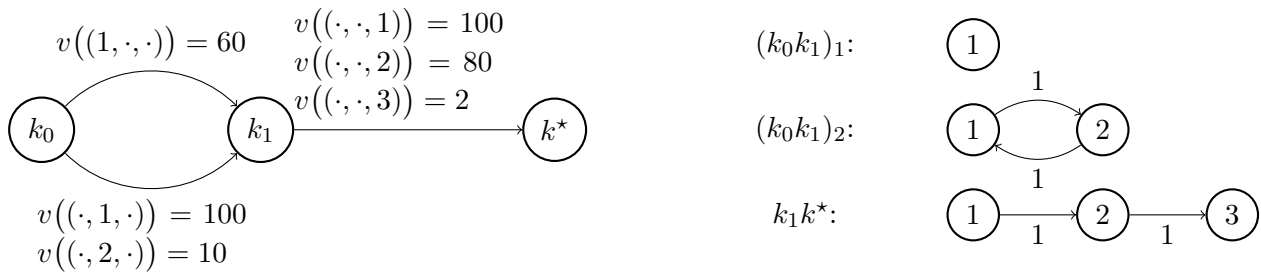

\begin{example}
\emph{
A vehicle wants to travel from $k_0$ to $k^\star$ in the network of Figure~\ref{fig:subpath},
in which each arc is 60 km long.
A state of the background process $B(t)$ is a 3-tuple consisting of the state 
of $(k_0k_1)_1$ (the upper link between $k_0$ and $k_1$), $(k_0k_1)_2$ 
(the lower link between $k_0$ and $k_1$)
and $k_1k^\star$, in that order.
A closer look at the velocities reveals that the speed on link $(k_0k_1)_1$ is always 60 km/h,
whereas the speeds on links $(k_0k_1)_2$ and $k_1k^\star$ can take 
two and three values, respectively.
For instance, the speed on link $(k_0k_1)_2$ is 100 km/h if this link is in state 1 and
10 km/h if this link is in state 2 (irrespective of the states of the other links).
We consider the situation that the initial state of the background process is set $B(t_0) = (1,1,1)$.
Using Theorem~\ref{Theorem1} we find the costs of using the upper and lower link:
\[
\mathbb{E}[\tau_{(k_0k_1)_1}^{B(t_0)}] = 1 \:\text{h}, \quad \quad \mathbb{E}[\tau_{(k_0k_1)_2}^{B(t_0)}] = 1.02\:\text{h}. 
\]
{\sc edsger} would therefore
advise to travel via $(k_0k_1)_1$. Comparison of the expected travel times of the
two paths between $k_0$ and $k^\star$ does however reveal that this is not optimal,
as these are 12.21h and 11.58h for traveling via $(k_0k_1)_1$ and 
$(k_0k_1)_2$ respectively.
Thus even though $(k_0k_1)_1$ is in the optimal path
to $k_1$, it is not the expected shortest path to $k^\star$.
A reason for this phenomenon lies in the fact that {\sc edsger} only takes the expectation
into account, while {variability} plays a role here as well, as this {variability}
directly relates to different conditions on future links. 
To see this, first note that traveling to $k_1$ via the upper link will always take 1 hour.
Upon arrival in $k_1$ there is a high probability that the background process of $k_1k^\star$ has
transitioned to a state with reduced speed.
If link $(k_0k_1)_2$ spends most time in state 1, traveling to $k_1$ via the 
lower link can be done within 1 hour.
In that case there is still a high probability that link $k_1k^\star$ can (partly) be traveled
at speed 100. Traveling via $(k_0k_1)_2$ therefore yields a lower expected travel time
than traveling via $(k_0k_1)_1$, despite the fact that we have $\mathbb{E}[\tau_{(k_0k_1)_1}^{B(t_0)}] < \mathbb{E}[\tau_{(k_0k_1)_2}^{B(t_0)}]$.
}\hfill$\Diamond$
\end{example}

A second reason for the potential loss in optimality lies in the the fact that
{\sc edsger} does not exploit its dynamical use in its fullest extent: due to the fact that in the DP approach one {\it knows} that one is allowed to change the path along the way, it potentially leads to lower cost than the dynamically used version of {\sc edsger}. This effect is highlighted in the example provided below.

\begin{example}
\emph{
The objective is to minimize the expected travel
time from $k_0$ to $k^\star$ in the network of Figure~\ref{fig:adaptinfo}
with $B(t_0) = (1,1,1,1)$.
Every arc has two states and a vehicle can drive 100 km/h on an arc
if the state of the arc is
1 and 80 km/h otherwise. 
The transitions from state 1 to state 2
have rate 0.1 and the transitions from state 2 to 
state 1 have rate 1 for all arcs. 
We first note that the expected travel
times on all paths are equal, such that
{\sc edsger} cannot
distinguish between these paths.
However, {\sc edsger} does not take into account that, in case $k_1$ is chosen as next node,
updated information about the speeds on the two optional paths from $k_1$ to $k^\star$
can be used to pick the optimal path given this updated information.
VI does use this information and as such picks the route via $k_1$,
yielding a 37 second improvement in expected travel time.
Note that, for a travel distance of 100 km, this is only a very small
improvement.
}\hfill$\Diamond$
\end{example}
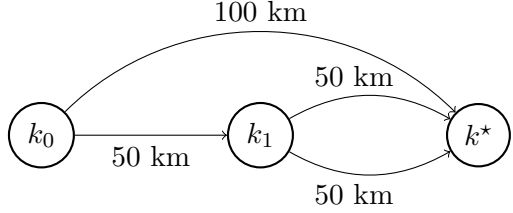
\begin{figure}[h]
\begin{center}
\begin{tikzpicture}[
roundnode/.style={circle, draw=black, thick, minimum size=7mm},
extranode/.style={circle, draw=black, thick, minimum size=4mm},
nonvisnode1/.style={circle, draw=white, thick, minimum size=1mm},
nonvisnode2/.style={circle, draw=white, thick, minimum size=17mm},
curnode/.style={circle, draw=red!60, fill=red!20, thick, minimum size=7mm},
squarednode/.style={rectangle, draw=red!60, fill=red!5, very thick, minimum size=5mm},
]
\node[roundnode]    (1)     {$k_0$};
\node[roundnode]    (2)     [right= 2 of 1]         {$k_1$};
\node[roundnode]    (3)     [right= 2 of 2]         {$k^\star$};

\draw[->, bend angle = 30, bend left] (2) to node [above,pos=0.4] {50 km} (3);
\draw[->, bend angle = 30, bend right] (2) to node [pos=0.4,below] 
{50 km} (3);
\draw[->, bend angle = 45, bend left] (1) to node [above,pos=0.5] {100 km} (3);
\draw[->] (1.east) -- (2.west) node [pos=0.5,below] 
{50 km};
\end{tikzpicture}
\caption{Example in which {\sc edsger} is not optimal,  
{as it does not fully exploit the possibility to
change the chosen route}.}
\label{fig:adaptinfo}
\end{center}
\end{figure}

Whereas the above instance shows that in principle we can construct cases in which DP leads to lower expected delay than {\sc edsger}, the experiments of Section~4 show that the difference is typically tiny (if there is a difference at all). Importantly, the modest loss of optimality that {\sc edsger} experiences is compensated by a 
potentially significant improvement in terms of the computational cost: the experiments show that the computational effort of {\sc edsger} can be substantially lower
than that of the VI algorithm.
This is because, in contrast to VI,
{\sc edsger} does not necessarily compute the expectation and
transition probabilities on all arcs in the network.



We now provide a detailed description of the shortest path algorithm that is used
within {\sc edsger}. 
As discussed above, {\sc edsger} uses this stochastic shortest path algorithm
to output a path from the current node to the destination.
The first node in this path defines the next action, i.e.,
the policy of {\sc edsger} is to travel to this node.

The shortest path 
algorithm is a stochastic version of the A$^\star $-algorithm.
{It} assigns a label to every node in the graph
and updates these labels iteratively. The initial labels are
$D_{k_0} = 0$ for the source $k_0$ and $D_k = \infty$ for $k \in N \setminus \{k_0\}$. 
A set ${\mathcal V}$ is used to store nodes with labels that are no longer
altered by {subsequent steps of the algorithm},
and is initialized by ${\mathcal V} := \emptyset$.
Every iteration of the algorithm
a new `current node' $c$ is chosen according to
\begin{align}
    c := \argmin_{k \in N \setminus {\mathcal V}}\{D_k + \mbox{\sc lb}_k\}. \label{eq:c}
\end{align}
Here $\mbox{\sc lb}_k$ is a lower bound
on the travel time from node $k$ to destination $k^\star$. 
Since the maximum speeds
and the distances of all arcs are known, 
lower bounds on the arc travel times can be computed in an elementary way. 
Applying Dijkstra's algorithm on a graph with
these lower bounds yields $\mbox{\sc lb}_k$ for $k \in N$.

With the current node $c$ chosen according to (\ref{eq:c}), the earlier
derived expressions (\ref{eq:solphi}) and (\ref{eq:trans})
are used to compute 
\begin{align}
     d_k &:= D_c + \sum_{s \in \mathcal{I}} p_c^s\,
    \mathbb{E}[\tau^s_{ck}] = D_c + {\bs p}_c^\top {\Phi}_{ck} \label{eq:dn} \\
    \label{L2} p_{k}^{s'} &:= \sum_{s \in \mathcal{I}} p_c^s\,
    \mathbb{P}(B(\tau_{ck}^s) = s')
\end{align}
for all $k \in \NB(c)$ and $s' \in \mathcal{I}$. 
In the expressions \eqref{eq:dn} and \eqref{L2} we let
$p_{k}^{s}$ 
denote the probability that upon arrival in $k \in N$
the background state is $s \in \mathcal{I}$, and
${\bs p}_{k} = (p_{k}^{s})_{s \in \mathcal{I}}$ is the vector of
these probabilities. In addition, $\Phi_{ck}$ is the $|{\mathcal I}|$-dimensional vector with entries $\mathbb{E}[\tau^s_{ck}]=\phi_s(d_{ck}\,|\,c,k)$ for $s\in{\mathcal I}.$
Note that, as $D_c$ can be seen as an estimate of the expected travel time from $k_0$ to
$c$, $d_k$ is an estimate for the expected travel time from $k_0$ to $k$.
The labels of the neighbors $k \in \NB(c) \setminus {\mathcal V}$ 
are now updated, i.e., if $d_k < D_k$ we set $D_k = d_k$ and store
the vector $p_k$ for this node (replacing a previous
stored value if present). {We also store the proposed path to $k$, 
which is the stored path to $c$ with $k$ itself appended.}
After performing this updating step, 
$c$ is added
to ${\mathcal V}$ and a new current node $c$ is 
chosen, again according to (\ref{eq:c}). 

The described
procedure is now repeated until
the destination node $k^\star$ is set as the current node $c$.
{The path the algorithm outputs is the stored path for $k^\star$.}
In the implementation of the algorithm 
a heap $H$ can be used to reduce the complexity of the algorithm, 
similar to the use of the heap in
Dijkstra's algorithm \cite{goldberg1996}. 
Pseudocode is provided in Algorithm~\ref{alg:edsger}.

\begin{algorithm}[h]
 \KwResult{path from $k_0$ to $k^\star$ given $B(0) = s$}
 Initialization: $D_{k_0} = 0, D_k = \infty$ for 
 $k \neq k_0$, heap $H = \{(D_{k_0} + \mbox{\sc lb}_{k_0} (\text{key}), p_{k_0}, k_0, \{k_0\} \text{ (path})\}, {\mathcal V} = \emptyset$\;
 \While{$H$ nonempty}{
 1. Extract tuple $(D_k + \mbox{\sc lb}_k,p_k,k, \text{path})$ with minimum first entry (key) from $H$\;
 2. \textbf{If} $k = k^\star$ quit and return path. \textbf{Else if} $k \in {\mathcal V}$ go back to step 1. \textbf{Else} continue\;
 3. \For{neighbor $k'$ in $\NB(k) \setminus {\mathcal V}$}{
 a. Compute $d_{k'}$ and $p_{k'}$ with (\ref{eq:solphi}) and (\ref{eq:trans})\;
 b. If $d_{k'} < D_{k'}$ set $D_{k'} = d_{k'}$
 and insert $(D_{k'} + \mbox{\sc lb}_{k'}, p_{k'}, k', \text{path } + \{k'\})$ in $H$\;
 }
 4. Add $k$ to ${\mathcal V}$
 }
 \caption{Outline shortest path algorithm in {\sc edsger}}
 \label{alg:edsger}
\end{algorithm}


Despite the fact that, in contrast to the VI
algorithm, {\sc edsger} does not necessarily compute the expectation
and transition probabilities for all arcs in the network,
the algorithm still suffers from the curse of dimensionality.
In this respect, note that the cost of evaluating the matrix exponential
in (\ref{eq:solphi})
will contribute substantially to the cost of {\sc edsger}. Moreover,
in case of a large state space, the dimension of $Q$ will make this
evaluation intractable. To overcome this drawback, we will introduce
the $\edsger$ algorithm. $\edsger$ can be regarded as an improved
version of {\sc edsger}, as it exploits the local-correlation
structure of the background process to reduce
the computational costs of {\sc edsger} drastically.

\subsection{$\edsger$ algorithm}
Both the VI algorithm and the {\sc edsger} algorithm do not use the
assumption of local-$r$-correlation (Assumption~1), which states that
the velocity on arc $k\ell \in A$ is only dependent on the congestion
levels of the arcs at most $r$ arcs away. The main idea behind the $\edsger$ algorithm
is to exploit this assumption, so as to overcome the curse of dimensionality.
Numerical examples in the next section will show that
$\edsger$ is indeed highly efficient and offers real-time response,
while still providing near-optimal results.

Similar to {\sc edsger}, at any decision epoch
{$\edsger$ uses a shortest path
algorithm that aims to}
output a minimal path from the current node to the
destination. 
The next arc in this path is traveled and
upon arrival in a new node
the shortest
path algorithm is called again.
We denote the resulting routing policy 
by $\pi_{\text{E}^\star}$.
The shortest path algorithm
used in $\edsger$ (Algorithm~2 below)
is very similar to the algorithm
used in {\sc edsger}, but
$\edsger$ exploits Assumption~1, in combination with a specific  approximation, 
to drastically speed up the computation of $d_{k'}$ and $p_{k'}$.

The main idea behind the $\edsger$-algorithm is that, 
under the assumption of local-$r$-correlation,
the velocity levels on an arc
are only dependent on the state of the 
background processes on the arcs in $A_{k\ell}^r$.
This implies that the expectation of the travel time
can be computed using just these processes.
As above, we write $A~=~\{k_1\ell_1, \dots, k_n\ell_n\}$.
Denote with $Q_{k\ell}^r$ the transition rate matrix of the 
process $B_{k\ell}^r(t) = (X_{k_i\ell_i}(t), k_i\ell_i \in A_{k\ell}^r)$ 
with state space $\mathcal{I}_{k\ell}^r$.

{The dynamics on arc $k\ell$ can be
completely  described by the process $B_{k\ell}^r(t)$,
as the velocity levels on arc $k\ell$ depend only on the
state of this process.
Thus, the velocity on arc $k\ell$ at time $t$, defined
as $v_{k\ell}((X_{k_i\ell_i}(t))_{k_i\ell_i \in A})$,
only depends on the values of entries $k_i\ell_i \in A_{k\ell}^r$.
We can therefore write the velocity on arc $k\ell$
at time $t$ as $v_{k\ell}((X_{k_i\ell_i}(t))_{k_i\ell_i \in A_{k\ell}^r})$.
Denote the truncation of $s = (s_{k_i\ell_i})_{k_i\ell_i \in A} \in \mathcal{I}$ to
the entries corresponding to links in $A_{k\ell}^r$ by 
$s_{k\ell}^r$, i.e., $s_{k\ell}^r = (s_{k_i\ell_i})_{k_i\ell_i \in A_{k\ell}^r}$.
Then we write $v_{k\ell}(s) = v_{k\ell}(s_{k\ell}^r)$ 
for the velocity level on arc $k\ell$
whenever $s \in \mathcal{I}$.}
\begin{theorem}
Denoting $\Phi^r(d \,|\, k,\ell)=(\mathbb{E}[\tau_{k\ell}^s(d)])_{s_{k\ell}^r
\in \mathcal{I}_{k\ell}^r}$, it holds that
\begin{align}
    \Phi^r(d \,|\, k,\ell)
    &= 
    \exp\left\{\begin{bmatrix}[1.3]
    d\,(V^r_{k\ell})^{-1}Q^r_{k\ell} & d\,(V^r_{k\ell})^{-1}{\boldsymbol 1} \\ {\boldsymbol 0}^\top & 0 \end{bmatrix}\right\}
    \begin{bmatrix}[1.3]
    {\boldsymbol 0} \\ 1
    \end{bmatrix}. \label{eq:phir}
\end{align}
Here $V^r_{k\ell}$ is a diagonal 
matrix with entries $v_{k\ell}(s_{k\ell}^r)$
for $s_{k\ell}^r \in \mathcal{I}_{k\ell}^r$.
\end{theorem}
\begin{proof}
The claim follows from Theorem~1, local-$r$-correlation (Assumption~1),
and the independence of the Markov processes
on the arcs.
\end{proof}
Particularly when $r$ is relatively small, this way of computing the expected per-edge travel times yields
significant computational savings, due to the fact that the dimension of $Q_{k\ell}^r$
no longer grows exponentially with the number of arcs.
This does however not directly yield a solution
to the curse of dimensionality, as the
computation of the transition probabilities still
involves the matrix $Q$.
We therefore propose to 
use the following (typically highly accurate) approximation for the transition
probabilities. 
Recall that the transition probabilities
of a Markov process $B(t)$ with transition rate matrix $A$ after a time $t \in \mathbb{R}$
can be expressed in terms of a matrix exponential:
\begin{align*}
    \mathbb{P}(B(t) = s \,| \,B(0) = s') = [e^{tA}]_{s',s}.
\end{align*}
{Note furthermore that we have a good estimate, $D_c$, for the travel time
from $k_0$ to $c \in N$}.
Based on {these two observations} we approximate $p_c^s$, the probability that $B_{k\ell}^r(t) = s$
upon arrival in $c$ given that $B(0) = s'$, by
\begin{align*}
    [e^{D_c Q_{ck}^r}]_{(s')_{k\ell}^r,s}.
\end{align*}
This results in the following expression for $d_k$:
\begin{align}
    d_k &:= D_c + 
    \sum_{s \in \mathcal{I}_{ck}^r} [e^{Q_{ck}^rD_c}]_{(s')_{k\ell}^r,s}\,
    \mathbb{E}\big[\tau^{s_{ck}^r}_{ck}\big] = D_c 
    + [e^{Q_{ck}^rD_c}]_{(s')_{k\ell}^r}
    \Phi^r(d_{ck}\,|\, c,k); \label{eq:edsgerstar}
\end{align}
cf.\ (\ref{eq:dn}). The corresponding pseudocode is given in Algorithm~2. 

\begin{algorithm}[h]
 \KwResult{path from $k_0$ to $k^\star$ given $B(0) = s$}
 Initialization: $D_{k_0} = 0, D_k = \infty$ 
 for $k \neq k_0$, heap $H = \{D_{k_0} + \mbox{\sc lb}_{k_0} \text{ (key)},k_0,
 \{k_0\} \text{ (path)}\}, {\mathcal V} = \emptyset$\;
 \While{$H$ nonempty}{
 1. Extract tuple $(D_k + \mbox{\sc lb}_k,k, \text{path})$ with minimum key from $H$\;
 2. \textbf{If} $k = k^\star$ quit and return path. \textbf{Else if} $k \in {\mathcal V}$ go back to step 1. \textbf{Else} continue\;
 3. \For{neighbor $k'$ in $\NB(k) \setminus {\mathcal V}$}{
 a. Compute $d_{k'}$ with (\ref{eq:edsgerstar})\;
 b. If $d_{k'} < D_{k'}$ set $D_{k'} = d_{k'}$
 and insert $(D_{k'} + \mbox{\sc lb}_{k'}, k', \text{path} + \{k'\})$ in $H$\;
 }
 4. Add $k$ to ${\mathcal V}$
 }
 \caption{Outline shortest path algorithm in $\edsger$}
 \label{alg:edsgerstar}
\end{algorithm}

\subsection{Implementation details}
\label{subsecimplementation}
We now discuss several implementation details for the VI-,
{\sc edsger}-, and $\edsger$-algorithms. 
{We will in particular show
how to rewrite specific computations in a form that allows the
application of efficient numerical functions. 
Moreover, we will recommend 
the use of efficient data structures}.
We start by discussing
the use of these techniques for the VI algorithm
and continue with the {\sc edsger}-, and the
$\edsger$-algorithms.


\paragraph{\it $\circ$~The value iteration algorithm.}
As discussed in Section~3.1,
this DP algorithm solves the Bellman optimality
equations and outputs a routing policy that minimizes our objective function.
It is an iterative procedure that
assigns values $\tilde{J}_0^{\pi^\star}(k,s)$
to all $(k,s) \in N \times \mathcal{I}$
and updates these values 
such that they
converge to $J^{\pi^\star}$, the set of optimal values.
Concretely, the algorithm sets
$\tilde{J}_0^{\pi^\star}(k^\star,s) = 0$
and $\tilde{J}_0^{\pi^\star}(k,s) = \infty$ for $k \in N, k \neq k^\star, s
\in \mathcal{I}$ and
iteratively updates these values according to the scheme 
\begin{align}
    \tilde{J}_i^{\pi^\star}(k,s) = \min_{\ell \in \NB(k)}\left\{\mathbb{E}[\tau_{k\ell}^s] + \sum_{s' \in \mathcal{I}}\,
    \mathbb{P}(B(\tau_{k\ell}^s) = s')\tilde{J}_{i-1}^{\pi^\star}(\ell,s')\right\},
    \label{eq:schemeVI}
\end{align}
cf.\ the scheme \eqref{eq:Vuiteen}.
Convergence of the iterative procedure has been widely studied in the literature; we refer to e.g.\
\cite{Bertsekas1991} and \cite{Bonet2007} and references therein. Specifically, it
has been proven that $\tilde{J}_i^{\pi^\star}(k,s)$ converges to $J^{\pi^\star}(k,s)$
for all $k \in N, s \in \mathcal{I}$.


The expectations and transition probabilities in the updating step \eqref{eq:schemeVI}
can be computed
by applying Theorem~\ref{Theorem1} and Corollary~\ref{cor:transprob}.
As was already 
noted in Remark~\ref{rem:suff},
the computation of the $(|\mathcal{I}| + 1) \times (|\mathcal{I}| + 1)$-dimensional
matrix exponential
\begin{align}
    \exp\left\{\begin{bmatrix}[1.3]
    d_{k\ell}V_{k\ell}^{-1}Q & d_{k\ell}V_{k\ell}^{-1}\boldsymbol{1} \\
    \boldsymbol{0}^\top & 0
    \end{bmatrix}\right\} \label{eq:matrixexp}
\end{align}
suffices.  
Note that the upper left $|\mathcal{I}| \times |\mathcal{I}|$-
and upper right 
$|\mathcal{I}| \times 1$-block
correspond to 
$\left(\mathbb{P}(B(X(\tau_{k\ell}^s) = s')\right)_{(s,s') \in 
\mathcal{I} \times \mathcal{I}}$ and
$\left(\mathbb{E}[\tau_{k\ell}^s]\right)_{s \in \mathcal{I}}$, respectively.
This means that, by writing 
$i_s$ for the index of $s$ in $\mathcal{I}$,
$\boldsymbol{p}_{k\ell}^s = (\mathbb{P}(B(\tau_{k\ell}^s) = s')_{s' \in \mathcal{I}}$
and
$\boldsymbol{\tilde{J}}_{i-1}^{\pi^\star}(\ell) = 
(\tilde{J}_{i-1}^{\pi^\star}(\ell, s'))_{s' \in \mathcal{I}}$,
the iterative step \eqref{eq:schemeVI} can be written as
\begin{align}
\tilde{J}_i^{\pi^\star}(k,s) &= \min_{\ell \in \NB(k)}\left\{
    \begin{bmatrix}
    (\boldsymbol{p}_{k\ell}^s)^\top
    & \mathbb{E}[\tau_{k\ell}^s]
    \end{bmatrix}
    \begin{bmatrix}[1.3]
    \boldsymbol{\tilde{J}}_{i-1}^{\pi^\star}(\ell) \\
    1
    \end{bmatrix}
    \right\}
    =
    \min_{\ell \in \NB(k)}\left\{
    \exp\left\{\begin{bmatrix}[1.3]
    d_{k\ell}V_{k\ell}^{-1}Q & d_{k\ell}V_{k\ell}^{-1}\boldsymbol{1} \\
    \boldsymbol{0}^\top & 0
    \end{bmatrix}\right\}_{i_s}
    \begin{bmatrix}[1.3]
    \boldsymbol{\tilde{J}}_{i-1}^{\pi^\star}(\ell) \\
    1
    \end{bmatrix}
    \right\}. \label{eq:schemeVIu}
\end{align}
Importantly, it now suffices to compute the {\it product of}\, the matrix exponential
\eqref{eq:matrixexp} and the vector $(\boldsymbol{\tilde{J}}_{i-1}^{\pi^\star}(\ell, k^\star), 1)^\top$
to evaluate the objective function in \eqref{eq:schemeVIu}. 
For such a product, also known as the \textit{action} of a matrix exponential,
most programming software include compiled functions,
examples of which are MatrixExp[\,\,] in Mathematica and
scipy.linalg.expm(\,\,) in Python; these functions are typically considerably faster than 
first computing  the matrix exponential
and the vector, and subsequently their product.

Besides the evaluation of the matrix exponential,
the computational costs of DP are strongly affected by
the data structure used for the $Q$-matrix.
Implementing the $Q$-matrix in a sparse way
significantly decreases the costs of constructing and storing this matrix.
Sparsity additionally reduces the costs of updating the values $\tilde{J}_i^{\pi^\star}(k,s)$,
since compiled matrix exponential functions in programming
software are generally faster in case the input is a sparse matrix.

\paragraph{\it $\circ$~The\, {\sc edsger}- and $\edsger$-algorithms.}
Also in the implementation of {\sc edsger} and $\edsger$,
a significant speed-up can be achieved by
(i)~working with the action of the matrix exponential and 
(ii)~exploiting the sparsity of the $Q$-matrix.
Here we will focus on the former speed-up, i.e., the one due to the use of the action of the matrix exponential, 
as the motivation for using a sparse $Q$-matrix is identical to the one
in the value iteration case.

Note that {\sc edsger} uses a matrix exponential in the computations
of $d_k$ and $\boldsymbol{p}_k$ (see \eqref{eq:dn} and \eqref{L2}). If we let
$\boldsymbol{p}_{k\ell}$ be the matrix 
$(\mathbb{P}(B(\tau_{k\ell}^s) = s')_{(s,s') \in \mathcal{I} \times 
\mathcal{I}}$, we can derive the following equivalence: 
\begin{align*}
    \begin{bmatrix} \boldsymbol{p}_{k}^\top & d_k-D_k \end{bmatrix}
    = 
    \boldsymbol{p}_c^\top \begin{bmatrix} \boldsymbol{p}_{k\ell} & \boldsymbol{\Phi}_{c\ell} \end{bmatrix}
    = 
    \begin{bmatrix} \boldsymbol{p}_c^\top & 0 \end{bmatrix}
    \begin{bmatrix}[1.3] \boldsymbol{p}_{k\ell} & \boldsymbol{\Phi}_{c\ell} \\ \boldsymbol{0}^\top & 1 \end{bmatrix}
    = 
    \begin{bmatrix} \boldsymbol{p}_c^\top & 0 \end{bmatrix}
    \exp\left\{\begin{bmatrix}[1.3]
    d_{k\ell}V_{k\ell}^{-1}Q & d_{k\ell}V_{k\ell}^{-1}\boldsymbol{1} \\
    \boldsymbol{0}^\top & 0
    \end{bmatrix}\right\}.
\end{align*}
Recall that for a vector $\boldsymbol{p}$ and square matrix $A$ for which $\boldsymbol{p}e^A$ exists, we have that 
$(\boldsymbol{p}e^A)^\top = (e^{A})^\top\boldsymbol{p}^\top = e^{A^\top}\boldsymbol{p}^\top$, yielding
\begin{align*}
    \begin{bmatrix}[1.3] \boldsymbol{p}_{k} \\ d_k \end{bmatrix}
    = \begin{bmatrix}[1.3] \boldsymbol{0} \\ D_k \end{bmatrix}
    + 
    \exp\left\{\begin{bmatrix}[1.3]
    d_{k\ell}Q^\top V_{k\ell}^{-1} & \boldsymbol{0} \\
    d_{k\ell}\boldsymbol{1}^\top V_{k\ell}^{-1} & 0
    \end{bmatrix}\right\}\begin{bmatrix}[1.3] \boldsymbol{p}_c \\ 0 \end{bmatrix}.
\end{align*}
{The same procedure can be followed 
to compute \eqref{eq:edsgerstar}
in the implementation of $\edsger$}. 


\subsection{Decreasing state space}
\label{subsecspeedup}
$\edsger$ yields significant computational
savings compared to DP (using value iteration) and the {\sc edsger} algorithm.
However, substantial reductions of the computational costs are possible. 
This is achieved by performing
preprocessing steps to reduce the network
size and the state space of the background process. These
are general speed-up techniques and can also be performed 
in the context of value iteration or {\sc edsger}. 

\subsubsection{Decreasing network size}
\begin{wrapfigure}{l}{0.35\textwidth}
  \begin{center}
    \includegraphics[width=\textwidth]{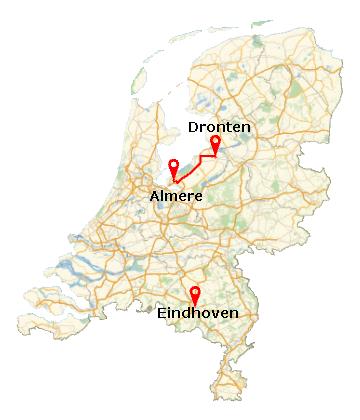}
  \end{center}
  \caption{Distance Almere-Dronten to Eindhoven}
  \label{fig:SmallRLargeN}
\end{wrapfigure}
Three specific ideas will be discussed in greater detail. 
The first one uses Yen's algorithm \cite{Yen1970, Yen1971}
to reduce the size of the network, by
deleting arcs that cannot be on shortest paths. 
The other
proposed ideas relate to decreasing the state space of the background 
process: 
the first idea considers hitting probabilities and excludes background
states for which the hitting probability is below a given threshold, whereas the second uses 
historical data to exclude background states which
are (extremely) rare events.

The first proposed speed-up technique decreases the network size by
deleting certain nodes and
arcs in the network.
The technique is motivated by the fact that not all
arcs in the network
will be considered by travelers.
The network 
in Figure~\ref{fig:SmallRLargeN}
depicts the shortest route (in distance) from
Almere to Dronten, in the Dutch road
system. In case there is congestion on this route
a traveler might wish to take a different route. 
One can argue, however, that the 
conditions in the network will 
never be such that the traveler
wishes to use an
arc around distant cities, such as, in our example, Eindhoven.
As a consequence, the roads around the
city of Eindhoven can be omitted when considering
the roads to travel from Almere to Dronten.

\begin{wrapfigure}{r}{.48\textwidth}
    \begin{minipage}{\linewidth}
    \centering\captionsetup[subfigure]{justification=centering}
    \begin{tikzpicture}[
roundnode/.style={circle, draw=black, thick, minimum size=7mm},
visnode/.style={circle, draw=black, thick, minimum size=7mm},
curnode/.style={circle, draw=red!60, fill=red!20, thick, minimum size=7mm},
squarednode/.style={rectangle, draw=red!60, fill=red!5, very thick, minimum size=5mm},
]
\node[roundnode]    (1)     {$k_0$};
\node[roundnode]    (2)     [right= 4 of 1]         {$k_1$};
\node[roundnode]    (3)     [right= 2 of 2]         {$k^\star$};

\draw[->] (1.east) -- (2.west) node [pos=0.5,above] {};
\draw[->] (2.east) -- (3.west) node [pos=0.5,above] {};
\draw[->] (1) to [out=30,in=150] (2);
\draw[->] (1) to [out=20,in=160] (2);
\draw[->] (1) to [out=10,in=170] (2);
\draw[->] (1) to [out=40,in=140] (2);
\draw[->] (1) to [out=-30,in=-150] (2);
\draw[->] (1) to [out=-20,in=-160] (2);
\draw[->] (1) to [out=-10,in=-170] (2);
\draw[->] (1) to [out=-40,in=-140] (2);
\draw[->] (1) to [out=-40,in=-160] (3);
\end{tikzpicture}
    \subcaption{Before Yen's algorithm}
    \label{fig:SSa}
    \begin{tikzpicture}[
roundnode/.style={circle, draw=black, thick, minimum size=7mm},
visnode/.style={circle, draw=black, thick, minimum size=7mm},
curnode/.style={circle, draw=red!60, fill=red!20, thick, minimum size=7mm},
squarednode/.style={rectangle, draw=red!60, fill=red!5, very thick, minimum size=5mm},
]
\node[roundnode]    (1)     {$k_0$};
\node[roundnode]    (2)     [right= 4 of 1]         {$k_1$};
\node[roundnode]    (3)     [right= 2 of 2]         {$k^\star$};

\draw[->] (1.east) -- (2.west) node [pos=0.5,above] {};
\draw[->] (2.east) -- (3.west) node [pos=0.5,above] {};
\draw[->] (1) to [out=20,in=160] (2);
\draw[->] (1) to [out=10,in=170] (2);
\draw[->] (1) to [out=-20,in=-160] (2);
\draw[->] (1) to [out=-10,in=-170] (2);
\end{tikzpicture}
    \subcaption{After Yen's algorithm}
    \label{fig:SSb}
    \includegraphics[height=5.1cm,keepaspectratio]{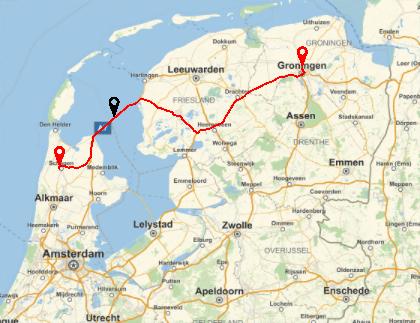}
    \subcaption{Real-world example}
    \label{fig:SSc}
\end{minipage}
\caption{Yen's algorithm 
deletes important alternative route}\label{fig:SS}
\end{wrapfigure}
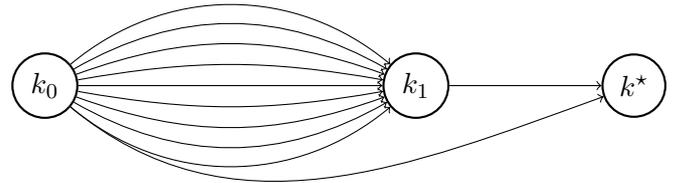
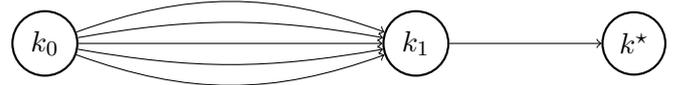
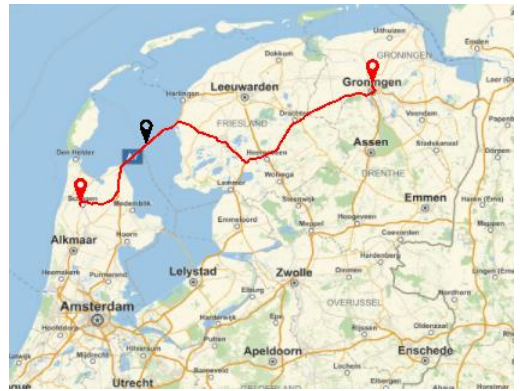
We now present a procedure, Yen's algorithm, that algorithmically determines the arcs that can be excluded.
The algorithm is first used to derive the $m$ shortest paths from source $k_0$ 
to destination $k^\star$ in a deterministic network, in which
a vehicle can drive at maximum velocity levels. 
Denote with $G' = (N',A')$ the network that solely consists of the
nodes and arcs in these $m$ shortest paths.
{Then we propose to reduce the network size by only considering
$G'$ in the routing problem.
Note that the value of $m$ should not be too large, as this
will not yield computational savings, but also not too small,
as this might eliminate arcs that are on the optimal route.
}

G\"uner et al.\ \cite{Guner2012} propose to reduce the network
$G$ to a network similar to $G'$, and compute a routing policy for this
reduced network.
We however propose an additional step in the construction of a reduced network, 
{in which we add edges to $G'$ which offer alternative routes 
in case these do not exist in $G'$.
An example is provided in} Figure~\ref{fig:SS}. 
Applying Yen's $m$-shortest
path algorithm with $m = 5$ in the graph of Figure~\ref{fig:SS}{\sc a}
yields Figure~\ref{fig:SS}{\sc b}, which shows that arc
$k_1 k^\star$ is present in every of the $5$ shortest paths
outputted by Yen's algorithm.
This is undesirable, since this arc
can be congested during the period
of travel, meaning a traveler might prefer arc $k_0k^\star$.
After an application of Yen's algorithm, $k_0 k^\star$ is however
no longer in the considered network and therefore not in
the routing policy. A real-world example is shown
in Figure~\ref{fig:SS}{\sc c}, where a vehicle wishes to travel
from Groningen to Schagen. Application of Yen's $m$-shortest path
algorithm may lead to $m$ paths that all contain the 
Afsluitdijk dam (black marker). 
A few times a year the Afsluitdijk
is closed for a few hours due to a major accident and 
in this case there are no alternative routes in $G'$.

{To avoid a situation as described above,
in which there is no alternative to a heavily congested arc, a second
step is introduced, in which arcs that
offer these alternative routes are
added to $G'$}.
To this end, it is first checked which arcs are in all of the
$m$ shortest paths. 
Yen's algorithm is then repeatedly
used to find $l$ shortest paths in the network in which one of these
arcs is deleted each time.
The sets of $l$ shortest paths are added to $G'$ 
{and this network
is used to construct a routing policy.}

To analyze the dynamics on $k\ell \in A'$, 
we only need information on the
state of the process $B_{k\ell}^r(t)$.
Decreasing the network size as described above
can therefore yield a significant reduction in
the size of the state space as well:
for $k\ell \in A \setminus (\cup_{i \in A'} A_i^r)$ 
the corresponding Markov processes $X_{k\ell}(t)$
give no information on the dynamics on $G'$
and can therefore be omitted in any further analysis.

\subsubsection{Use of bounds on hitting probabilities}\label{DHP}
One technique for reducing
the state space of the background process is 
based on hitting probabilities. The idea is to only
include, for a small number $\epsilon$, background states $s \in \mathcal{I}$
for which $\mathbb{P}(B(t) = s) > \epsilon$ for some $t \leqslant T$,
with $T >0$ denoting the time of arrival
at destination $k^\star$.
Since this time horizon is evidently not known in advance,
one could work with a $M > 0$ that serves as a crude upper bound on $T$.
Given an initial
background state $s \in \mathcal{I}$ an upper bound for
the hitting probabilities is then given by
\begin{align}
    \mathbb{P}(\exists t \in [0,M]: B(t) = s' \,|\, B(0) = s) \leqslant
    \prod_{i = 1}^{|A|} \mathbb{P}(\exists t \in [0,M]:
    X_{k_i\ell_i}(t) = s_{k_i\ell_i}' \,|\, X_{k_i\ell_i}(0) = s_{k_i\ell_i}) \label{eq:ub}
\end{align}
In case $\mathcal{S}_{k_i\ell_i} = \{1,2\}$ for some $i \in \{1,\dots,|A|\}$
with transition rates $\lambda_i,\mu_i$, we have
\begin{align*}
    f_{12}(M) := \mathbb{P}(\exists t \in [0,M]:
    X_{k_i\ell_i}(t) = 2 \,|\,  X_{k_i\ell_i}(t)(0) = 1)
    \leqslant \mathbb{P}(Y_{\lambda_i}(t) \leqslant M) = 1-e^{-\lambda_i M},
\end{align*}
where $Y_{\lambda_i}(t) \sim \text{Exp}(\lambda_i)$. We equivalently have that
$f_{21}(M) \leqslant 1-e^{-\mu_i M}$,
and by definition $f_{11}(M) = f_{22}(M) = 1$.

Whenever $\mathcal{S}_{k_i\ell_i} = \{1,\dots,n_{k_i\ell_i}\}$ with
$n_{k_i\ell_i} > 2$, we can rely on standard results on sums of independent exponentially distributed random variables \cite{Bibinger2013} in order to find
an upper bound on the probabilities in (\ref{eq:ub}). An example
is given in Appendix \ref{appendix:c}, in which it is shown how
to use these results in case the Markov process on an arc is a birth-death process.
Now, if (an upper bound of) the derived bound from (\ref{eq:ub})
is smaller than some predefined $\epsilon > 0$, we do not
consider the corresponding background state $s' \in \mathcal{I}$
in our further analysis.

\subsubsection{{Use of historical data}}\label{DHD}
A technique to further 
reduce the background state space  
is directly based on  historical data of the process $B(t)$ on $G$.
When this data is available, it can be analyzed
to determine which states $s \in \mathcal{I}$ are most common
to occur in reality and which states have never occurred 
during the time period during which the historical data was recorded.
For example, in a network with 50 arcs in which each
arc has two states (congested and uncongested), the
situation in which {\it all} arcs are congested can
theoretically occur, but will not be observed
in real road networks, as will  be confirmed
by the historical data. 
Hence a strategy could be to eliminate 
all states $s \in \mathcal{I}$ that have never been attained
by the process $B(t)$, so as to reduce 
the state space of this
background process.

The three techniques discussed in this section can be used as preprocessing steps in
dynamic routing. In case of {\sc edsger} and 
$\edsger$, they can be applied prior to  
every call of the shortest path algorithm. 
An outline of the resulting procedure is given in Algorithm~\ref{alg:total}.

\begin{algorithm}[H]
 \KwResult{Travel policy from $k_0$ to $k^\star$}
 initialization: network $G = (N,A)$, $B(0) = s$, Node = $k_0$\;
 Step 1: Preliminary node deletion\;
 a. Derive $m$ shortest paths on $G$ to form $G'$\;
 b. Identify arcs $\{k_1\ell_1,\dots,k_n\ell_n\}$ that are on all shortest paths\;
 c. Derive $l$ shortest paths on $G \setminus k_i\ell_i$ for $i = 1,\dots,n$
 and add to $G'$\;
 d. Reduce state space background process by looking at deleted arcs\;
 Step 2: Preliminary background state deletion\;
 a. Delete states with hitting probability constraint\;
 b. Delete states by historical data analysis\;
 Step 3: Use shortest path algorithm\;
 a. Determine path from $k_0$ to $k^\star$\;
 b. Travel to outputted next node. Stop if this is $k^\star$.
 Otherwise return to Step 1 with state upon
 arrival as initial state\;
 \caption{Outline Dynamic Routing with Preprocessing steps}
 \label{alg:total}
\end{algorithm}

\section{Numerical Experiments}
\label{sec:numexp}
This section presents a series of numerical experiments that 
demonstrate that the $\edsger$ algorithm
performs near-optimally with high efficiency.
That is, corresponding to the earlier
introduced notions of 
{\it distance-to-optimaliy} and {\it efficiency},
$\edsger$ outputs a value
close to the minimally achievable value
while essentially being real-time.
To substantiate this claim, we have considered a broad range of traffic scenarios and networks of various dimensions.

\begin{wrapfigure}{r}{0.48\textwidth}
  \begin{center}
    \includegraphics[width=\textwidth]{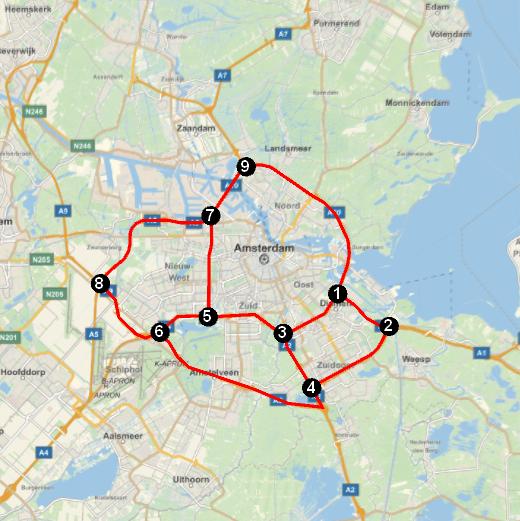}
  \end{center}
    \caption{Amsterdam highway system}
    \label{fig:AdamHS}
\end{wrapfigure}
More specifically, we first consider a small network 
and show that value iteration (VI), {\sc edsger}, and $\edsger$
outperform two deterministic algorithms (in terms of distance-to-optimality),
in case of a simple as well as a more sophisticated background process.
Second, we 
increase the size of the network to show that the run-time of
VI and {\sc edsger} grows exponentially in the network size, whereas
the run-time of $\edsger$ is substantially less affected.
Importantly, $\edsger$ still yields close-to-optimal results in these larger
networks. Last, we consider routing in a network of realistic size
and show that $\edsger$ is still highly efficient and nearly optimal.


The first two experiments 
consider routing on
the highway system around Amsterdam.
They
compare the distance-to-optimality 
of VI, {\sc edsger} and $\edsger$
with two deterministic algorithms.
The two deterministic algorithms, used as benchmarks, 
both employ the A$^\star$-algorithm.
The first deterministic algorithm, 
as before referred to as `Deterministic Static' (abbreviated to DS),
applies the A$^\star$-algorithm on the network with
{\it maximum} velocities. 
As the maximum speeds are fixed,
the A$^\star$-algorithm is executed once to determine the complete
travel path, which is followed until the 
destination is reached. The second deterministic algorithm,
which we will call `Deterministic Dynamic' (abbreviated to DD),
does use the available information on the velocities in the network. 
Every intersection the algorithm calls the A$^\star$-algorithm on 
a network with the {\it current} velocities to determine
the next arc to travel. 

Experiment~1 uses a Markovian background process in which
every arc has just two states: congested or uncongested.
It is shown that VI, {\sc edsger} and $\edsger$ 
outperform the deterministic algorithms in terms of 
distance-to-optimality in this simple setup.
Experiment~2 considers a more detailed
background process in which the Markov process of each
arc can attain three possible states (so as to more accurately 
model the speeds the vehicles can drive at). 
The background process in addition includes
a global event with Erlang distributed holding time
(which could represent a reasonably predictable change in the global 
circumstances, e.g.\ rush hour; recall Remark~\ref{rem:erlang}).
Graphical and numerical summaries show that 
the difference in distance-to-optimality
between our algorithms
and the deterministic algorithms is even more significant 
than in Experiment~1.
Comparing efficiency shows that the deterministic algorithms
and $\edsger$ can be executed in real-time, 
contrary to VI and {\sc edsger}, whose
computational costs suffer from the size of the state space. 

In Experiment~3 we evaluate the  computational costs
of the various routing algorithms as function of the network size. We do so by working with 
an elementary network model, that is then extended with additional arcs to assess its impact on the algorithms' speeds. In addition, we quantify the effect of the dimension of the background process.
The final experiment, Experiment~4,
considers routing on the entire Dutch highway network
and  highlights the influence of various model parameters on the
performance of the algorithms. Moreover, the experiment considers
the speed-up techniques introduced in Section~\ref{subsecspeedup}.
For the experiments we implemented the networks and
routing algorithms in Wolfram Mathematica 12.0 on 
an Intel\textregistered\ Core\texttrademark\ i7-8665U 1.90GHz computer.

\subsection*{Experiment 1}
To get a first impression of the performance of the various algorithms, we
start by considering a relatively small network. 
Our findings reveal that in this network, using  a background process in which every link
has just two states, both {\sc edsger} and $\edsger$ efficiently yield
close-to-optimal results 
Consider the network of the Amsterdam
highway system (Figure~\ref{fig:AdamHS}), with 9 intersections and 24 links
between these intersections (i.e., 12 bidirectional arcs).
We pick the following framework: 
\begin{itemize}
    \item[$\circ$] The background process of every arc contains just two states, uncongested (corresponding to state $1$) and
    congested (corresponding to state $2$). Transition rates are tuned with NDW data.
    \item[$\circ$] The state of an arc only affects
    the speeds on the directly attached arcs, i.e., there is local-$1$-correlation.
    \item[$\circ$] There is no process $Y$ that induces global correlation.
    \item[$\circ$] In case there is no incident on the arc we let the vehicle speed be 100 km/h if  there is no incident on the directly adjacent arcs, and 80 km/h otherwise.
    In case there is an incident on the arc we let the vehicle speed be 40 km/h if 
    there is no incident on the directly adjacent arcs, and 20 km/h otherwise.
    \item[$\circ$] In the network there is a maximum of three incidents simultaneously, to 
    bound the size of the state space and guarantee tractability of VI.
\end{itemize}

Consider a traveler interested in minimizing the total
expected travel time from node 1 to node 8.
We compare the DP policy with the 
routing policies under {\sc edsger}, 
$\edsger$ and the two deterministic algorithms DS and DD.
The DP policy is derived from the VI-algorithm,
using the implementation guidelines as described in Section~\ref{subsecimplementation}.
Implementation of the two deterministic policies follows 
from a standard implementation of the A$^\star$-algorithm.
The policies of {\sc edsger} and $\edsger$
can be found by storing
the output of their shortest path algorithms
for every initial
state $(k_0,s) \in N \times \mathcal{I}$. 
The example below, intended to demonstrate the principles 
underlying {\sc edsger} in a concrete setting, shows that the policy 
of {\sc edsger} in node 1, with as background state that every
arc except for the arc from node 2 to node 1 is uncongested,
is to travel to node 3.

{\it Example.}
{
We use {\sc edsger} to route from
node 1 to node 8 in Figure~\ref{fig:AdamHS},
when upon leaving the only congested arc is
the arc from node 2 to node 1.
As a first step, Dijkstra is used on a network with maximum speeds,
to determine the lower bounds $\mbox{\sc lb}_k$ of the travel time from
node $k \in N$ to node 8, which in this case are given by
\begin{align*}
    \mbox{\sc lb}_1 = 0.18, \enspace \mbox{\sc lb}_2 = 0.23, \enspace \mbox{\sc lb}_3 = 0.13, \enspace \mbox{\sc lb}_4 = 0.16, \enspace \mbox{\sc lb}_5 = 0.09,
    \enspace \mbox{\sc lb}_6 = 0.05, \enspace \mbox{\sc lb}_7 = 0.09, \enspace \mbox{\sc lb}_8 = 0, \enspace \mbox{\sc lb}_9 = 0.13
\end{align*}
Now the shortest path algorithm within {\sc edsger} is used
to determine the next node to travel to.
Node 1 is set as current node and 
we initialize $D_1 = 0, D_2 = \dots = D_9 = \infty$.
Then $\exp(d_{1k'}V^{-1}Q)$ is computed to determine 
$d_{k'}$ and $\boldsymbol{p}_{k'}$ in \eqref{eq:dn} and \eqref{L2} for $k' = 2,3,9$, the neighbors of 1, to give:
\begin{align*}
    d_2 = 0.19, \enspace d_3 = 0.04, \enspace d_9 = 0.10.
\end{align*}
Since $d_{k'} < D_{k'}$ for $k' = 2,3,9$ we set
$D_{k'} = d_{k'}$ and store ($D_{k'} + \mbox{\sc lb}_{k'}, p_{k'}, k', \{1,k'\}$)
for all neighbors $k'$ of 1. The updated labels now yield:

\begin{table}[h]
\begin{tabular}{|c|c|c|c|c|c|c|c|c|}
\hline
     $D_1 = 0$ & $D_2 = 0.19$ & $D_3 = 0.04$ & $D_4 = \infty$ & $D_5 = \infty$ & 
     $D_6 = \infty$ & $D_7 = \infty$ & $D_8 = \infty$ & $D_9 = 0.10$ \\ \hline
\end{tabular}
\end{table}
The new current node is then $\argmin_{k' \in N \setminus \{1\}}\{D_{k'} + \mbox{\sc lb}_{k'}\} = 3$.
Computing $d_{k'}$ and $p_{k'}$ for $k' = 4,5$ gives
\begin{align*}
    d_4 = 0.07, \enspace d_5 = 0.08.
\end{align*}
We therefore store ($d_{k'} + \mbox{\sc lb}_{k'}, p_{k'}, k', \{1,3,k'\}$)
for $k' = 4,5$ and update the labels:

\begin{table}[h]
\begin{tabular}{|c|c|c|c|c|c|c|c|c|}
\hline
     $D_1 = 0$ & $D_2 = 0.19$ & $D_3 = 0.04$ & $D_4 = 0.07$ & $D_5 = 0.08$ & 
     $D_6 = \infty$ & $D_7 = \infty$ & $D_8 = \infty$ & $D_9 = 0.10$ \\ \hline
\end{tabular}
\end{table}
The next current node is set as
$\argmin_{k' \in N \setminus \{1,3\}}\{D_{k'} + \mbox{\sc lb}_{k'}\} = 5$.
This procedure is repeated until the current node is set as 8.
The paths stored for node 8 is $(1,3,5,6,8)$, 
indicating the first node to travel to is node 3.
}
\hfill$\Diamond$

\begin{figure}[h]
    \centering
    \includegraphics[width = 0.5\textwidth]{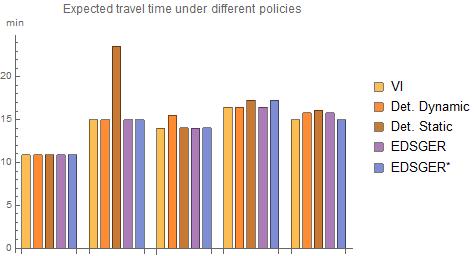}
    \caption{Highlight a few background states}
    \label{fig:exp1highlight}
\end{figure}

\begin{table}[h]
    \centering
    \begin{tabular}{|l||c|c|c|c|c||r|}\hline
    & Average & 1 Inc. & 2 Inc. & 3 Inc. & WA & Run-time \\
    Value iteration & 12.75 & 11.60 & 12.26 & 12.83 & 11.16 & 2.032 \\
    {\sc edsger} & 12.76 & 11.60 & 12.27 & 12.84 & 11.16 & 0.094 \\
    $\edsger$ & 12.76 & 11.61 & 12.27 & 12.84 & 11.16 & 0.023 \\
    Det. Dynamic & 12.83 & 11.64 & 12.33 & 12.91 & 11.17 & 0.003 \\
    Det. Static & 13.08 & 11.61 & 12.41 & 13.19 & 11.17 & $<$ 0.001 \\
    \hline
    \end{tabular}
    \caption{Costs algorithms}
    \label{tab:exp1comp}
\end{table}

Figure~\ref{fig:exp1highlight} shows the expected travel time (in minutes)
of the derived policies in five initial background states. 
The first set of five bars corresponds to a state of complete non-congestion, i.e.,
all of the arcs are uncongested upon departure. 
The expected travel times under the different policies
are in this case similar.
Note that this is not surprising, 
as the influence of unpredictable events is minimal:
the distance between
node 1 and node 8 is small and therefore 
the probability of occurrence of an incident
on the shortest path as well.
The second set of bars shows the largest 
difference in expected travel time 
between VI and DS.
This difference in expected travel time
is significant, the expectation under the deterministic
policy being more than 1.5 as large as the expectation under VI.
The third, fourth and fifth set of bars show 
the largest difference
between the expected travel time under VI
and the expected travel time under 
DD, $\edsger$ and {\sc edsger} respectively. 
Note that these differences are small,
and their policies are thus close-to-optimal. 

The fact that DD, {\sc edsger} 
and $\edsger$ yield close-to-optimal 
results in the framework of this experiment can also
be seen in Table~\ref{tab:exp1comp}:
\begin{itemize}
    \item[$\circ$] The first column  shows the expected travel time
from node 1 to node 8 under the different policies, averaged (evenly) over
all possible initial background states.
 \item[$\circ$] The second column
 shows the expected travel time from node 1 to node 8, averaged (evenly)
 over all initial background states in which one incident has occurred.
 Columns three and four show the same for respectively two and three incidents.
 \item[$\circ$] The fifth column provides a weighted average (WA) of the expectations corresponding to 
all possible initial background states; for a given initial background state,
the weight equals its limiting probability.
 \item[$\circ$]
The last column contains the run-time of the algorithms (in sec.).
\end{itemize}
The average and weighted averages of {\sc edsger} and $\edsger$
are close to the results of VI. 
DD performs relatively well, whereas the 
values under DS are noticeably suboptimal.
Comparison in run-time shows that the computational costs of VI are
an order larger than those of the other algorithms.
In Experiment~3 we will investigate the
efficiency of the algorithms in greater detail, and show
that the difference in computational costs
becomes more substantial when increasing the
maximum number of incidents in the network.

We conclude that in this simple framework, in which every arc
can only have two states, {\sc edsger} and $\edsger$ 
yield nearly optimal results, 
while being roughly one order of magnitude faster than VI.
It was also shown that {\sc edsger} and $\edsger$ 
perform better than
the two deterministic algorithms DS and DD.
Especially in case there are incidents in the network,
the difference in distance-to-optimality is significant.
As the occurrence of incidents is relatively rare, the 
distance-to-optimality of the four algorithms
is similar if there are no incidents in the network upon departure.
The next experiment shows that the difference-to-optimality
typically grows when adding more detail to the model.

\subsection*{Experiment 2}
In this experiment we consider a more involved
background process, and show that {\sc edsger}
and $\edsger$ still yield nearly optimal results. The two algorithms
outperform the deterministic algorithms (in terms
of distance-to-optimality), and are at the same time considerably faster than 
VI.
The example, moreover, demonstrates the comprehensiveness
of our model, by illustrating the possibility of (i) adding recurrent
events and (ii) working with more velocities per arc.

We again consider the network of the Amsterdam
highway system (Figure~\ref{fig:AdamHS}), but extend the 
background process of Experiment~1.
Concretely, we have three speeds per arc (rather than two), and include
a global event.
We include the extra speed level and global event
in the following way:
\begin{itemize}
    \item[$\circ$] Every arc has three states, uncongested (state $1$), 
    congested (state $2$) and recovery (state $3$). A link can transition from
    an uncongested to a congested state, but not vice versa: 
    from a congested
    state a link must first enter the recovery state before it
    returns to the uncongested state. 
    The recovery state represents the time between the clearance
    of an incident and the time the traffic conditions return
    to the free-flow speed.
    \item[$\circ$] The state of an arc only affects
    the speeds on the attached arcs, i.e., there is local-$1$-correlation.
    \item[$\circ$] There is a recurrent event
    that induces global correlation, to be interpreted as a rush hour.
    This event affects the speeds on the roads on the inner circle of the
    highway system, i.e., all arcs between nodes 1, 3, 5, 7 and 9.
    We will denote these arcs as category I arcs, whereas we will
    refer to the other arcs as category II.
    As a rush hour is a recurring event the time till its onset has a relatively low variance. That is the reason why we chose to not model this time by an 
    exponential distribution but rather by an 
    Erlang distribution  (as pointed out in Remark~\ref{rem:erlang}). In our experiments we took
    four phases. Only in the last phase, which we identify as the start of the
    rush hour, the speeds on the arcs in category I are affected.
    \item[$\circ$] 
    Speed levels for the different scenarios can be found in Table~\ref{tab:speedsexp2}.
    \item[$\circ$] In the network there is again a maximum of three non-uncongested links simultaneously, to 
    bound the size of the state space and guarantee tractability of VI.
\end{itemize}

\begin{table}[h]
    \centering
    \begin{tabular}{|*6c|}
\hline
& &  \multicolumn{2}{c}{Non-rush} & \multicolumn{2}{c|}{Rush}\\
\hline
& & 0   & $\geqslant 1$    & 0   & $\geqslant 1$ \\
\multirow{3}{*}{Cat. I} & 1   &  100 & 80   & 70 & 60\\
& 2   &  40 & 20   & 20  & 10 \\
& 3   &  70 & 50 & 50 & 40\\
\hline
\multirow{3}{*}{Cat. II} & 1   &  100 & 80   & 100 & 80\\
& 2   &  40 & 20 & 40 & 20 \\
& 3   &  70 & 50 & 70 & 50 \\
\hline
\end{tabular}
    \caption{Speed levels on the arcs in km/h. The rows specify 
    the category and state of the arc, the columns specify
    if there is rush hour or not, and distinguish between the numbers of attached
    arcs that is not in state 1.
    }
    \label{tab:speedsexp2}
\end{table}

\begin{figure}[h]
    \centering
\begin{tikzpicture}[shorten >=1pt,node distance=2cm,on grid,auto]
    \node[state] (q_0) {$1$};
    \node[state] (q_1) [right=of q_0] {$2$};
    \node[state] (q_3) [right=of q_1] {$3$};
    \node[state] (q_4) [right=of q_3] {$4$};

    \path[->]
    (q_0) edge [] node {$\lambda$} (q_1)
    (q_1) edge [] node {$\lambda$} (q_3)
    (q_3) edge [] node {$\lambda$} (q_4);
\end{tikzpicture}
    \caption{The four Erlang phases of Rush hour}
    \label{fig:erlang}
\end{figure}
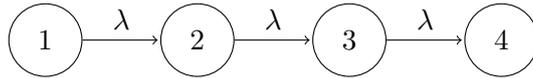



\begin{figure}[h]
    \centering
    \includegraphics[width = 0.5\textwidth]{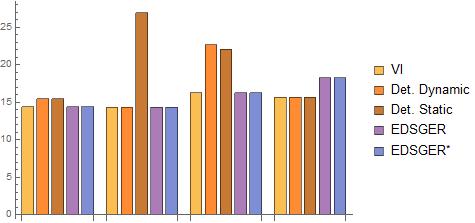}
    \caption{Highlight a few background states}
    \label{fig:exp2highlight}
\end{figure}

Figure~\ref{fig:exp2highlight} shows the expected travel time (in minutes) of the algorithms
in four initial background states. The first set of five bars again corresponds 
to a state of complete non-congestion. 
The  expected  travel  times  under  the  different  policies  are  in  
this  case  similar. The second and third set of 
bars shows the largest difference in expected travel time 
between VI and DS and DD, respectively.
This difference in expected travel time is significant, especially for DD. 
The fourth set of bars show the largest difference between
the expected travel time under VI and the 
expected travel time under both $\edsger$ 
and  {\sc edsger}.  Note that this difference is 
small, and hence the policies {\sc edsger} and $\edsger$ are
close to optimal.

\begin{table}[h]
    \centering
    \begin{tabular}{|l||c|c|c|c|c|c|c|c||r|}\hline
    & \multicolumn{2}{c}{Phase 1} & \multicolumn{2}{c}{Phase 2} 
    & \multicolumn{2}{c}{Phase 3} & \multicolumn{2}{c||}{Phase 4} & Run-time \\ \hline
    & A & WA & A & WA & A & WA 
    & A & WA & \\
    Value iteration & 15.55 & 14.56 & 15.70 & 14.59 & 15.75 & 14.60 & 15.77 & 14.61 & 41.513 \\
    {\sc edsger} & 15.68 & 14.59 & 15.73 & 14.59 & 15.76 & 14.61 & 15.78 & 14.61 & 13.154\\
    $\edsger$ & 15.70 & 14.60 & 15.74 & 14.60 & 15.77 & 14.61 & 15.78 & 14.61 & 0.154 \\
    Det. Dynamic & 16.20 & 15.53 & 16.76 & 16.33 & 17.20 & 17.00 & 15.87 & 14.61 & 0.003 \\
    Det. Static & 17.64 & 15.73 & 18.56 & 16.59 & 19.27 & 17.30 & 19.57 & 17.62 & $<$0.001\\ \hline
    \end{tabular}
    \caption{Average (A) and weighted average (WA) expected travel time of the algorithms in the four
    different initial phases of the rush hour. Run-time (in sec.) reported as well.}
    \label{tab:exp2comp}
\end{table}

The low distance-to-optimality of {\sc edsger} and $\edsger$ 
can also be observed
from Table~\ref{tab:exp2comp}.
The table shows the average and weighted average
expected travel time (in minutes) under the different policies
and in the four different initial Erlang phases of
the event rush hour. Similar to Experiment~1, the weights
are set equal to the limiting probabilities.
The average and weighted averages of {\sc edsger} and $\edsger$
are very close to the results of VI. 
DS also performs well, whereas the 
values obtained by DD are noticeably suboptimal.
The last column of Table~\ref{tab:exp2comp} shows the run-time (in seconds)
of the different algorithms. Note that the computational costs of VI and
{\sc edsger} are substantially larger than those of $\edsger$
and the two deterministic algorithms.

Thus, summarizing the results of the first two experiments,
$\edsger$ has the best performance
in terms of distance-to-optimality and run-time.
More precisely, $\edsger$ outperforms the deterministic algorithms
in terms of distance-to-optimality,
and outperforms VI and {\sc edsger} in terms of
efficiency. The next example will show that the savings in run-time
will become even more pronounced when increasing the network size and/or the  maximum number of incidents. 

\subsection*{Experiment 3}
\begin{figure}[h]
\centering
     \begin{subfigure}[b]{0.15\textwidth}
     \centering
    \begin{tikzpicture}[
roundnode/.style={circle, draw=black, thick, minimum size=7mm},
visnode/.style={circle, draw=black, thick, minimum size=7mm},
curnode/.style={circle, draw=red!60, fill=red!20, thick, minimum size=7mm},
squarednode/.style={rectangle, draw=red!60, fill=red!5, very thick, minimum size=5mm},
]
\node[roundnode]    (1)     {S};
\node[roundnode]    (2)     [right= 1 of 1]         {};
\node[roundnode]    (4)     [below= 1 of 1]         {};
\node[roundnode]    (5)     [right= 1 of 4]         {D};

\draw[-] (1.east) -- (2.west) node [pos=0.5,above] {};
\draw[-] (4.east) -- (5.west) node [pos=0.5,above] {};
\draw[-] (1.south) -- (4.north) node [pos=0.5,above] {};
\draw[-] (2.south) -- (5.north) node [pos=0.5,above] {};
\end{tikzpicture}
    \caption{$2\times2$-network}
    \label{fig:Toymodel2}
\end{subfigure}
     \begin{subfigure}[b]{0.30\textwidth}
     \centering
    \begin{tikzpicture}[
roundnode/.style={circle, draw=black, thick, minimum size=7mm},
visnode/.style={circle, draw=black, thick, minimum size=7mm},
curnode/.style={circle, draw=red!60, fill=red!20, thick, minimum size=7mm},
squarednode/.style={rectangle, draw=red!60, fill=red!5, very thick, minimum size=5mm},
]
\node[roundnode]    (1)     {S};
\node[roundnode]    (2)     [right= 1 of 1]         {};
\node[roundnode]    (3)     [right= 1 of 2]         {};
\node[roundnode]    (4)     [below= 1 of 1]         {};
\node[roundnode]    (5)     [right= 1 of 4]         {};
\node[roundnode]    (6)     [right= 1 of 5]         {};
\node[roundnode]    (7)     [below= 1 of 4]         {};
\node[roundnode]    (8)     [right= 1 of 7]         {};
\node[roundnode]    (9)     [right= 1 of 8]         {D};

\draw[-] (1.east) -- (2.west) node [pos=0.5,above] {};
\draw[-] (2.east) -- (3.west) node [pos=0.5,above] {};
\draw[-] (4.east) -- (5.west) node [pos=0.5,above] {};
\draw[-] (5.east) -- (6.west) node [pos=0.5,above] {};
\draw[-] (7.east) -- (8.west) node [pos=0.5,above] {};
\draw[-] (8.east) -- (9.west) node [pos=0.5,above] {};
\draw[-] (1.south) -- (4.north) node [pos=0.5,above] {};
\draw[-] (4.south) -- (7.north) node [pos=0.5,above] {};
\draw[-] (2.south) -- (5.north) node [pos=0.5,above] {};
\draw[-] (5.south) -- (8.north) node [pos=0.5,above] {};
\draw[-] (3.south) -- (6.north) node [pos=0.5,above] {};
\draw[-] (6.south) -- (9.north) node [pos=0.5,above] {};
\end{tikzpicture}
    \caption{$3\times3$-network}
    \label{fig:Toymodel3}
\end{subfigure}
\caption{Examples of networks of size $n \times n$. Here `S' is the source, and `D' the destination.}
\label{fig:elementarynetworks}
\end{figure}
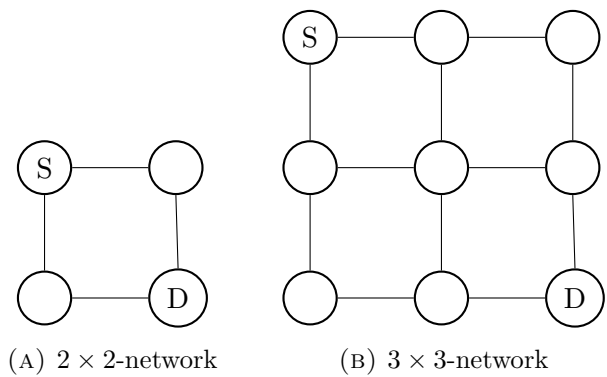

We consider elementary networks to 
assess the sensitivity of the run-time as a function of the network size,
and in addition as a function of the maximum number of incidents.
In the previous examples, in which the network size
and maximum number of incidents were small, we already noted that the 
computational costs of both VI
and {\sc edsger} are considerably higher than those of $\edsger$.
To directly demonstrate the computational savings when using $\edsger$,
we use networks of the type depicted in Figures~\ref{fig:elementarynetworks}{\sc a} and
\ref{fig:elementarynetworks}{\sc b}.
These networks show a $2 \times 2$- and $3 \times 3$-structure, but extending the model in the obvious way produces an $n \times n$-structure for any $n \geqslant 2$. 
We increase the value of $n$, so as to evaluate the impact on the computational costs for VI, {\sc edsger}, 
and $\edsger$.
For simplicity, the conditions on the network are set 
identical to those of Experiment~1. This includes the condition that there can
only be three incidents in the network simultaneously. However, below we 
also assess how increasing this bound on the number
of incidents affects the run-time of the algorithms.

\begin{figure}[h]
    \centering
    \begin{subfigure}[b]{0.40\textwidth}
 \centering
 \includegraphics[width=\textwidth]{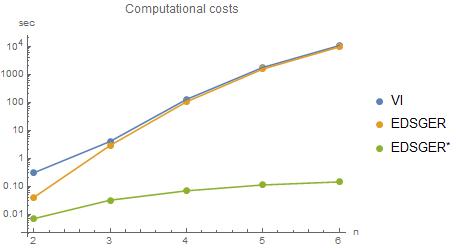}
 \caption{Computational costs in $n \times n$-graphs}
 \label{fig:nbyn}
\end{subfigure}
\begin{subfigure}[b]{0.40\textwidth}
 \centering
 \includegraphics[width=\textwidth]{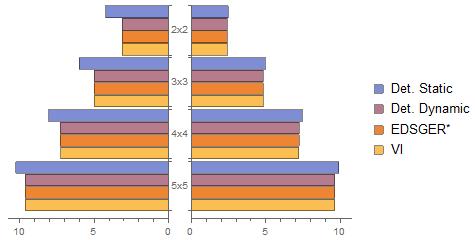}
 \caption{Performance algorithms in $n \times n$-graphs}
 \label{fig:nbynperformance}
\end{subfigure}
\caption{Effect of increasing the network size.}
\label{fig:increasenetworksize}
\end{figure}

Figure~\ref{fig:increasenetworksize}{\sc a} shows that, contrary to the
run-time of $\edsger$,
the run-times of both VI and {\sc edsger}
grow exponentially with the network size (note the
logarithmic scale).
This exponential increase 
can be explained by the fact that the size of the background
process, and thus the size of $Q$, grows exponentially with the
size of the network.
Since $\edsger$ only
uses the part of the state space that corresponds to the states on the
directly attached links, its computational costs 
are hardly affected by the network size.
We conclude that in larger networks
VI and {\sc edsger} 
become intractable, whereas $\edsger$ 
still offers real-time response. 

\begin{figure}[h]
    \centering
\begin{subfigure}[b]{0.40\textwidth}
 \centering
 \includegraphics[width=\textwidth]{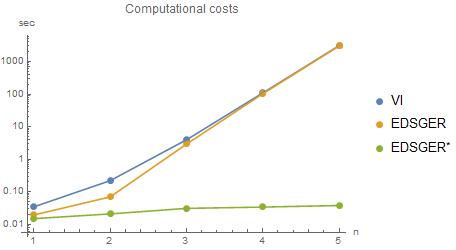}
 \caption{Computational costs with
 max. of $n$ incidents}
 \label{fig:3by3inc}
\end{subfigure}
\begin{subfigure}[b]{0.40\textwidth}
 \centering
 \includegraphics[width=\textwidth]{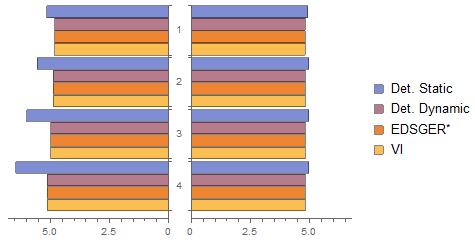}
 \caption{Performance algorithms in case
 $\leqslant n$ incidents}
 \label{fig:3by3incperf}
\end{subfigure}
\caption{Effect of increasing the number of incidents.\label{fig:incincrease}}
\end{figure}

Importantly, the substantial reduction of the run-time when using $\edsger$
(instead of VI, that is) does not correspond to a significant increase of the objective function:
Figure~\ref{fig:increasenetworksize}{\sc b} shows that $\edsger$ yields close to optimal
results. 
The figure shows the average
and weighted average expected travel time (in minutes), with weights
again chosen as the stationary probabilities, for different
network sizes. Observe that DD also yields
close-to-optimal results in this framework. This can, however,
be explained by the fact that the considered
instance is relatively simple; adding more detail to the framework,
as we did in Experiment~2, again leads to a more pronounced suboptimality
of this algorithm.

\begin{wrapfigure}{r}{0.4\textwidth}
  \includegraphics[width=\textwidth]{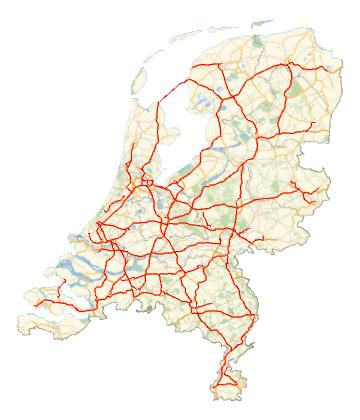}
     \caption{Dutch highway system}
    \label{fig:DutchHS}
\end{wrapfigure}
We observe the same behavior (in terms of computational costs
and the value of the objective function) when, instead of the network size,
the bound on the
maximum number of incidents is increased (see Figure~\ref{fig:incincrease}). 
That is, considering the 
$3 \times 3$-network displayed in Figure~\ref{fig:elementarynetworks}{\sc b},
the run-time of VI and {\sc edsger} grows exponentially with
the maximum number of incidents. Moreover, the run-time of
$\edsger$ is hardly affected by the increase of the maximum number
of incidents, while still performing nearly optimally.
The substantial 
difference in run-time is again due to the rapid growth of the dimension of $Q$
as function of the number of incidents. In contrast, the matrices
$Q^r_{k\ell}$, as used by $\edsger$, involve significantly fewer arcs,
and are therefore less affected by the number of incidents.
Note that, in case of local-1-correlation, the computational
costs of $\edsger$ will e.g.\ no longer increase if the number
of incidents exceeds the maximum degree of the graph.

From Figures~\ref{fig:increasenetworksize} and \ref{fig:incincrease} we conclude that
VI and {\sc edsger} are inefficient, making these
algorithms unsuitable for practical purposes. The results
of $\edsger$ are considerably more promising, as these showed 
that its run-time is hardly affected by the network size or maximum number of incidents,
while performing close-to-optimal in terms of cost.
This is why in our last experiment we investigate the tractability of $\edsger$ in
networks of realistic size. The experiment also considers the influence of
the model parameters and demonstrates the use of the speed-up techniques
introduced in Section~\ref{subsecspeedup}.


\subsection*{Experiment 4}
To confirm the feasibility of $\edsger$ 
in large networks,
the final experiment considers routing in a network
of realistic size: the Dutch highway
network (Figure~\ref{fig:DutchHS}),
with 93 intersections and 262 directed roads between these
intersections.
Besides confirming the feasibility of $\edsger$,
the experiment has the following three objectives:
\begin{itemize}
\item[$\circ$] Show that $\edsger$
still yields accurate results
in terms of expected travel time, supporting
our conclusions of the experiments above.
We will omit VI and {\sc edsger} in our analysis,
as these algorithms are intractable in large networks,
as demonstrated in Experiment~3.
Instead, we use simulations to compare the
travel time under $\edsger$, DD and DS with the travel
time under an optimal path;
\item[$\circ$]
Assess the effect of the parameters on the distance-to-optimality
of $\edsger$. Concretely, we compare the results of
$\edsger$ with the results of DD and DS for different parameter values;
\item[$\circ$]
Demonstrate how the speed-up techniques
described in Section~\ref{subsecspeedup} can be used in this network, and how
the different techniques affect the distance-to-optimality and the performance of
the algorithms.
\end{itemize}

To show the travel time under $\edsger$ is indeed close-to-optimal, 
we simulate the travel time for three OD-pairs.
The considered OD-pairs, denoted by `Short', `Medium'
and `Long', differ notably in length, with distances
22.5, 98.5 and 191.7 km respectively.
For simplicity we again use the same background process
as in Experiment~1.
We simulate realizations of this background process and
determine for every realization the travel time of $\edsger$, 
as well as the travel times of DS and DD.
These travel times are compared with the optimal travel time, i.e.,
the minimal achievable travel time under the given realization
of the background process. 
To measure the difference in travel time, we compute the
\textit{average percentage loss in travel time}:
$$
    \text{Average}\Big(\begin{frac}
    {\text{travel time algorithm} - \text{optimal travel time}}{\text{optimal travel time}}
    \end{frac} \cdot 100\% \Big).
$$
Note that a high value of this measure can also
arise for an algorithm that is optimal in terms
of expected travel time, e.g.\ VI.
This is due to the difference between
(i)~the notion of being (close-to-)optimal in expectation 
and~(ii)~the notion of being optimal for a given realization. 
Using this measure we can therefore only compare algorithms,
not assess individual values.

Table~\ref{tab:exp4obj1} shows the results of the simulations.
We first note that $\edsger$ is indeed feasible, as the displayed run-times $t$ 
(in sec.) are sufficiently low. The columns denoted by
U$\%$ and W$\%$ show the average percentage loss in travel time 
for a given OD-pair,
with initial states chosen uniformly and weighed (with the corresponding limiting probabilities),
respectively.

\begin{table}[h]
    \centering
    \begin{tabular}{|l||c|c|r|c|c|r|c|c|r|}\hline
    & \multicolumn{3}{c|}{Short} & \multicolumn{3}{c|}{Medium} & \multicolumn{3}{c|}{Long} 
    \\ \cline{2-10}
    & U$\%$ & W$\%$ & \multicolumn{1}{c|}{$t$} & U$\%$ & W$\%$ & \multicolumn{1}{c|}{$t$} 
    & U$\%$ & W$\%$ & \multicolumn{1}{c|}{$t$} 
    \\ \hline
    $\edsger$ 
    & 6.5 & 0.4 & 0.09 & 5.8 & 0.4 & 0.14 & 3.4 & 0.3 & 0.22 \\
    DD 
    & 7.3 & 0.4 & 0.01 & 5.9 & 0.8 & 0.02 & 4.9 & 0.7 & 0.03 \\ 
    DS
    & 25.7 & 1.3 & $<$0.01 & 5.1 & 0.5 & $<$0.01 & 3.0 & 0.4 & $<$0.01 \\
    \hline
    \end{tabular}
    \caption{Average percentage loss in travel time (W$\%$)
    and run-time of the algorithms ($t$); background process as in Experiment~1.}
    \label{tab:exp4obj1}
\end{table}

Table~\ref{tab:exp4obj1} shows that $\edsger$ is close-to-optimal
in this large network. If the initial states are drawn
according to the limiting probability (W$\%$), the average
percentage loss in travel time
is below $1\%$. In case the
initial state is chosen uniformly, the drawn initial states contain typically
many congested links. Even in this more extreme setting, the 
average percentage loss in travel time
is below $7\%$. 
Note that the results of DD are relatively close to the results of $\edsger$,
which is not surprising, as we have chosen the same framework as in Experiment~1 
(in terms of the structure of the background process).
Table~\ref{tab:exp4obj2} illustrates that 
the difference in distance-to-optimality between $\edsger$ 
and DD becomes more significant if we add more detail
to the background process (by using a setup as in Experiment~2).
In Tables~\ref{tab:exp4obj1} and \ref{tab:exp4obj2} 
the results of DS show
the disadvantage of not taking into account
any background information, 
mostly notable in case of the first OD-pair.

\begin{table}[h]
    \centering
    \begin{tabular}{|l||c|r|}\hline
    & \multicolumn{2}{c|}{Short} \\ \cline{2-3}
    & W$\%$ & \multicolumn{1}{c|}{$t$} \\ \hline
    $\edsger$ 
    & 2.0 & 0.60 \\
    DD 
    & 8.5 & 0.01 \\ 
    DS
    & 11.8 & $<0.01$ \\
    \hline
    \end{tabular}
    \caption{Average percentage loss in travel time (W$\%$) 
    and run-time of the algorithms ($t$); background process as in Experiment~2.}
    \label{tab:exp4obj2}
\end{table}

We  can in addition assess the effect of the parameters on the distance-to-optimality 
of the algorithms.
To this end, we again pick the setup of Experiment~1 and simulate, 
for different values of the incident rate $\alpha$ 
and clearance rate $\beta$, the travel time under $\edsger$, DD
and DS for the OD-pair `Short'.
Initial states are chosen according to their corresponding limiting probabilities.
Table~\ref{tab:exp4obj3}
shows the average percentage loss in travel time
for $\edsger$, DD and DS, respectively.

\begin{table}[h]
    \begin{tabular}{|c|c|c|c|}
    \multicolumn{4}{c}{$\edsger$} \\
    \hline
    \diagbox{$\alpha$}{$\beta$} & $10^{\minus 4}$ & 1 & 100 \\ \hline
    $10^{\minus 4}$ & 0.00 & 0.00 & 0.00 \\
    1 & 0.00 & 10.75 & 0.05 \\
    100 & 0.00 & 0.22 & 1.50 \\ \hline
    \end{tabular}
    \quad
    \begin{tabular}{|c|c|c|c|}
     \multicolumn{4}{c}{DD} \\
    \hline
    \diagbox{$\alpha$}{$\beta$} & $10^{\minus 4}$ & 1 & 100 \\ \hline
    $10^{\minus 4}$ & 0.00 & 0.00 & 0.00 \\
    1 & 0.00 & 12.58 & 1.30 \\
    100 & 0.00 & 1.46 & 29.28 \\ \hline
    \end{tabular}
    \quad
    \begin{tabular}{|c|c|c|c|}
     \multicolumn{4}{c}{DS} \\
    \hline
    \diagbox{$\alpha$}{$\beta$} & $10^{\minus 4}$ & 1 & 100 \\ \hline
    $10^{\minus 4}$ & 39.83 & 0.03 & 0.00 \\
    1 & 0.00 & 38.52 & 0.05 \\
    100 & 0.00 & 0.22 & 1.50 \\ \hline
    \end{tabular}
    \caption{Average percentage loss for $\edsger$ (left),
    DD (middle) and DS (right) and different values of the incident rate $\alpha$ and clearance rate $\beta$.}
    \label{tab:exp4obj3}
\end{table}

The tables reflect the non-stochastic nature of DS and the fact that
DD does not take any information on the duration of incidents into account.
A few observations we can do:
\begin{itemize}
    \item[$\circ$] $\edsger$ and DD perform extremely well in case the time
    between accidents is long (i.e., the incident rate is low). The probability
    of an accident occurring while traveling is in this case very low, thus
    any route circumventing existing incidents will suffice for a low travel time;
    \item[$\circ$] A small clearance rate (i.e., it will take very long
    before an accident has been cleared) has a potentially  
    high impact on the travel times of DS.
    In case $\alpha = 1$ or $\alpha = 100$, the limiting probability
    for congestion is very high, resulting in many congested arcs.
    Because there are so little uncongested arcs in the network,
    the optimal travel time is with high probability attained on the shortest path (in km),
    explaining the small average
    percentage loss of DS;
    However, if $\alpha = 10^{\minus 4}$, in the initial state approximately
    half of the arcs will be congested. DS does not reroute in case of
    incidents, and thus has a very high travel time in case an accident occurs on
    the chosen path;
    \item[$\circ$] DD is clearly suboptimal in case
    the incident and clearance rate take higher values. 
    This can be explained by the fact that DD 
    will reroute in case an incident has occurred on 
    a favorable path, whereas the high clearance rate
    will imply the incident duration to be short and
    thus the impact of the incident on the shortest path minimal.
    \item[$\circ$] 
    The highest value in the table of $\edsger$ arises in 
    case $\alpha = \beta = 1$.
    Note that this loss is smaller than for DD and DS,
    implying $\edsger$ still performs better than the deterministic
    algorithms. The relatively high loss in travel time
    for $\edsger$ is, as argued above, due to the difference between
    the notion of being (close-to-)optimal in expectation 
    and the notion of being optimal for a given realization.
    The loss under VI would likely be
    of the same order as that of $\edsger$.
    Computing the loss of both algorithms in the network
    of Experiment~1 with $\alpha = \beta = 1$ results
    e.g.\ in losses $7.00\%$ (VI) and $7.27\%$ ($\edsger$). 
\end{itemize}

As stated above, the third objective of this experiment is to demonstrate
the effect of the speed-up techniques described in Section~\ref{subsecspeedup} 
on the distance-to-optimality and the efficiency of the algorithms.
The outcome of this experiment, in which
we study the same routes as above (i.e., `Short', `Medium', `Long'), 
is shown in Table~\ref{tab:exp4speedup}.
We use Yen's $n$ shortest path algorithm to
decrease the network size and distinguish in Table~\ref{tab:exp4speedup}
between the different values of $n$, namely $n = 1,5,10,25$.
As argued in Section~\ref{subsecspeedup}, this decrease in network
size directly implies a decrease in the state space, as we do no
longer track the states of the arcs not contained in the reduced network.
The derived bound on the hitting probabilities in \eqref{eq:ub} 
gives that, for the OD-pair of `Short',
the probability of occurrence of four accidents
while driving has probability smaller
than $\epsilon = 0.0005$. 
Background states that contain four additional incidents with respect
to the initial state are therefore deleted.
A similar argument can be used to reduce the size of the state space
for the OD-pairs of `Medium' and `Long'.

\begin{table}[h]
    \centering
    \begin{tabular}{|l||c|c|c|c|c|c||c|c|c|c|c|c|}\hline
    \multirow{2}{*}{} & \multicolumn{6}{c||}{$n = 1$} & \multicolumn{6}{c|}{$n = 5$} 
    \\
    \cline{2-13}
    & \multicolumn{2}{c|}{Short} & \multicolumn{2}{c|}{Medium} & \multicolumn{2}{c||}{Long} 
    & \multicolumn{2}{c|}{Short} & \multicolumn{2}{c|}{Medium} & \multicolumn{2}{c|}{Long}
     \\ \hline
     nodes 
     & \multicolumn{2}{c|}{8} & \multicolumn{2}{c|}{18} & \multicolumn{2}{c||}{35} 
     & \multicolumn{2}{c|}{9} & \multicolumn{2}{c|}{18} & \multicolumn{2}{c|}{35} 
     \\ 
     arcs 
     & \multicolumn{2}{c|}{9} & \multicolumn{2}{c|}{22} & \multicolumn{2}{c||}{42} 
     & \multicolumn{2}{c|}{12} & \multicolumn{2}{c|}{23} & \multicolumn{2}{c|}{43} 
         \\ \hline
     & $\%$ & $t$ & $\%$ & $t$ & $\%$ & $t$ 
     & $\%$ & $t$ & $\%$ & $t$ & $\%$ & $t$ 
      \\ \hline
    $\edsger$ 
    & 0.8 & 0.04 & 0.5 & 0.07 & 0.5 & 0.12
    & 0.8 & 0.04 & 0.5 & 0.07 & 0.4 & 0.19
   \\
    DD 
    & 0.8 & 0.01 & 0.9 & 0.01 & 0.6 & 0.04
    & 0.8 & 0.01 & 0.9 & 0.01 & 0.5 & 0.10
    \\
    DS 
    & 1.9 & $<$0.01 & 0.6 & $<$0.01 & 0.6 & $<$0.01
    & 1.9 & $<$0.01 & 0.6 & $<$0.01 & 0.5 & $<$0.01 
    \\ \hline
    & \multicolumn{6}{c||}{$n = 10$}
    & \multicolumn{6}{c|}{$n = 25$} 
    \\   \cline{2-13}
    & \multicolumn{2}{c|}{Short} & \multicolumn{2}{c|}{Medium} & \multicolumn{2}{c||}{Long}
    & \multicolumn{2}{c|}{Short} & \multicolumn{2}{c|}{Medium} & \multicolumn{2}{c|}{Long}
    \\ \hline
    nodes 
    & \multicolumn{2}{c|}{10} & \multicolumn{2}{c|}{19} & \multicolumn{2}{c||}{37} 
     & \multicolumn{2}{c|}{16} & \multicolumn{2}{c|}{22} & \multicolumn{2}{c|}{40} 
     \\
     arcs
     & \multicolumn{2}{c|}{15} & \multicolumn{2}{c|}{25} & \multicolumn{2}{c||}{48} 
     & \multicolumn{2}{c|}{27} & \multicolumn{2}{c|}{35} & \multicolumn{2}{c|}{55} 
     \\
     \hline
     & $\%$ & $t$ & $\%$ & $t$ & $\%$ & $t$
     & $\%$ & $t$ & $\%$ & $t$ & $\%$ & $t$
     \\
     \hline 
     $\edsger$
      & 0.8 & 0.05 & 0.5 & 0.07 & 0.4 & 0.24
    & 0.8 & 0.05 & 0.5 & 0.09 & 0.4 & 0.34 \\
    DD
    & 0.8 & 0.01 & 0.9 & 0.01 & 0.5 & 0.14
    & 0.8 & 0.01 & 0.9 & 0.02 & 0.5 & 0.25 \\ 
    DS 
    & 1.9 & $<$0.01 & 0.6 & $<$0.01 & 0.5 & $<$0.01
    & 1.9 & $<$0.01 & 0.6 & $<$0.01 & 0.5 & $<$0.01\\ \hline
    \end{tabular}
    \caption{Average percentage loss ($\%$) and run-time ($t$) for
    different OD-pairs and different values of $n$.}
    \label{tab:exp4speedup}
\end{table}

The table reveals that the number of nodes
and arcs in the network increases with $n$
and consequently the run-time $t$ of the algorithms as well.
A moderate value of $n$ is therefore preferred.
This is also justified by the distance-to-optimality 
of the algorithms for moderate values of $n$.
That is, the table shows that for smaller values of $n$
the algorithms still yields close-to-optimal results, 
with the average percentage loss below $2\%$
for all cases. 
Note that the travel times for $n = 10$ and $n = 25$
are even similar to the case of $n = 5$.
This again advocates the choice of a small $n$, e.g.\ between
$n = 1$ and $n = 5$, 
as this reduces the run-time while 
the distance-to-optimality of the algorithms is not affected.

\section{Concluding Remarks}
In this paper we developed a new mechanism for describing the evolution of the velocities in a road traffic network, capable of modeling both recurrent and non-recurrent congestion. For this flexible velocity model we developed a routing algorithm that aims at minimizing the expected travel time. Extensive experiments showed that it outperforms competing models, in that it provides near-optimal results but at the same time offers real-time response.

Regarding the velocity model, we advocate the use of a continuous-time Markovian background process
to describe the speed driven by vehicles on the arcs in the network. We have developed various ways to deal with the potentially high dimension of this process. 
Future research could concern further operationalizing this model, with a specific focus on tuning the model's parameters using measurement data.
Concretely, we have argued that the proposed background process is able to incorporate the
influence of random events, as well as the influence of (nearly)
deterministic patterns, on the vehicle speed.
Insight into the occurrence of these events,
and their effect on the velocities,
can be provided relying on Intelligent
Transportation Systems (ITS).
However, it s not a priori evident how to map the information provided by a typical ITS on the parameters of our Markov model.
A future study could explore such calibration issues.

Regarding the  routing algorithms, we  have evaluated 
these on the basis of (i) distance-to-optimality (i.e., minimizing the expected travel time
between source and destination) 
and (ii) efficiency (i.e., run-time).
Numerical experiments justify the advice to use
$\edsger$ as routing algorithm, as $\edsger$ 
is close-to-optimal and has real-time run-time.
In the implementation of $\edsger$, as well as the
implementation of the other presented algorithms, 
the guidelines described in Section~\ref{subsecimplementation}
were used.
There may be opportunities to further 
speed-up the algorithm. 
For instance, parallel computing could potentially 
be used to speed up VI, {\sc edsger}, and $\edsger$.
Since the costs of computing the per-arc expected 
travel times are under $\edsger$ significantly lower
than under {\sc edsger} and VI, also when applying  parallel computing
$\edsger$ will outperform the competing 
algorithms.

Numerical experiments were conducted to show that
$\edsger$ efficiently yields close-to-optimal results 
under a broad range of realistic traffic scenarios.
The model parameters have been calibrated from
historical data on vehicle speeds and flows, 
collected by the National Data Warehouse for
Traffic Information (NDW) in the Netherlands.
Future research {will} focus on developing 
a more formal statistical 
procedure to estimate these parameters,
performing a detailed analysis 
of the occurrence and consequences of incidents;
cf.\ the earlier study by Snelder \cite{Snelder2013}.

\section*{Acknowledgements}
The authors would like to thank dr.\ Maaike Snelder (TU Delft)
for her helpful feedback on the manuscript as well as her suggestions
on the analysis of the NDW data.

\bibliographystyle{plain}
\bibliography{refs.bib}

\appendix
\section{Proof of Proposition 1}
\label{appendix:a}
\begin{repproposition}{prop1}
    Let $k\ell \in A$ and denote with $\tau_{k\ell}^{t}(d)$ the travel time on arc
    $k\ell$ for traveling a distance $d \in [0,d_{k\ell}]$ when leaving
    at $t \geqslant 0$. Then 
    $t_1 \leqslant t_2$ implies that $t_1 + \tau_{k\ell}^{t_1}(d)
    \leqslant t_2 + \tau_{k\ell}^{t_2}(d)$.
\end{repproposition}
\begin{proof}
If $\int_{t_1}^{t_2} v_{k\ell}( X(v))\,{\rm d}v \geqslant d$
it follows that
$$
t_1 + \tau_{k\ell}^{t_1}(d) = t_1 + 
\min\left\{t\geqslant 0: 
\int_{t_1}^{t_1+t} v_{k\ell}( X(v))\,{\rm d}v
\geqslant d \right\}
\leqslant t_1 + (t_2 - t_1) = t_2 \leqslant t_2 + \tau_{k\ell}^{t_2}(d).
$$
If $\int_{t_1}^{t_2} v_{k\ell}( X(v))\,{\rm d}v < d$
it follows that
\begin{align*}
t_1 + \tau_{k\ell}^{t_1}(d) &\leqslant t_1 + \min\left\{t\geqslant 0: 
\int_{t_2}^{t_1 + t} v_{k\ell}( X(v))\,{\rm d}v
\geqslant d\right\} 
\leqslant t_1 + 
\min\left\{t + t_2 - t_1 \geqslant 0: 
\int_{t_2}^{t_2 + t} v_{k\ell}( X(v))\,{\rm d}v
\geqslant d\right\} \\
&= t_1 + (t_2 - t_1) + \tau_{k\ell}^{t_2}(d) = t_2 + \tau_{k\ell}^{t_2}(d).
\end{align*}
This completes the proof.
\end{proof}

\section{Proof of Theorem 1}
\label{appendix:b}
\begin{reptheorem}{Theorem1}
Given a graph $G = (N,A)$ with a pair of nodes $k,\ell \in N$ and a distance $d\geqslant 0$,
it holds that
\begin{align*}
    V_{k\ell}\, \Phi'(d \,|\, k,\ell) &= 1 + Q\,\Phi(d \,|\, k,\ell), \\
    V_{k\ell}\, \Psi'(d \,|\, \alpha,k,\ell) &= (Q + \alpha I)\,\Psi(d \,|\, \alpha,k,\ell),
\end{align*}
with $V_{k\ell} := {\rm diag}\,\{(v_{k\ell}(s))_{s \in \mathcal{I}}\}$
and ${\boldsymbol 1}$ a $|\mathcal{I}|$-dimensional column vector of ones.
A solution for this system of linear differential equations can be written as
\begin{align}
    \Phi(d \,|\, k,\ell)
    &= 
    \begin{bmatrix} {\boldsymbol 1}^\top & 0 \end{bmatrix}
    \exp\left\{\begin{bmatrix}dV_{k\ell}^{-1}Q & dV_{k\ell}^{-1}{\boldsymbol 1} \\ {\boldsymbol 0}^\top & 0 \end{bmatrix}\right\}
    \begin{bmatrix}
    {\boldsymbol 0} \\ 1
    \end{bmatrix}, \\
    \Psi(d \,|\, \alpha,k,\ell) &=
    \exp\{d\,V_{k\ell}^{-1}(Q + \alpha I)\},
\end{align}
with ${\boldsymbol 0}$ an $|\mathcal{I}|$-dimensional
column vector of zeroes.
\end{reptheorem}
\begin{proof}
The proof uses a type of `infinitesimal argumentation' that is frequently relied upon in the context of fluid storage systems. 
Let $d \in [0,d_{k\ell}]$, $s,s' \in \mathcal{I}$.
Conditioning on a possible jump of the background process 
in $[0,\Delta]$, as $\Delta\downarrow 0$, recalling that scenarios with more than one jump have a probability that is $o(\Delta)$, 
\begin{align*}
    \phi_{s}(d\,|\,k,\ell) 
    &= (1 + \Delta q_{ss}){\mathbb E}\big[\tau^s_{k\ell}(d) \,|\,
    B(t) = s \enspace \forall t \in [0,\Delta] 
    \big] \\
    \quad &\hspace{5mm}+ \sum_{\tilde{s} \neq s} \Delta q_{s\tilde{s}}\,
    {\mathbb E}\big[\tau^s_{k\ell}(d) \,|\, \exists t \in [0,\Delta]: B(u) = s \enspace
    \forall u \in [0,t), B(u) = s' \enspace \forall u \in [t,\Delta)
    \big] + o(\Delta) \\
    &= \Delta + (1 + \Delta q_{ss})\phi_{s}(d - v_{k\ell}(s)\Delta \,|\,k,\ell)
    \\
    \quad &\hspace{5mm}+ \sum_{s' \neq s} \frac{\Delta q_{ss'}}{1 - e^{-q_{ss'}\Delta}} \int_{0}^\Delta 
    \phi_{s'}(d - v_{k\ell}(s)t - v_{k\ell}(s')(\Delta - t) \,|\,k,\ell)
    \,q_{ss'}e^{-q_{ss'}t}\,{\rm d}t + o(\Delta),
\end{align*}
where it is used that, for $s' \neq s$,
\begin{align*}
    {\mathbb E}\big[&\tau^s_{k\ell}(d) \,|\, \exists t \in [0,\Delta]: B(u) = s \enspace
    \forall u \in [0,t), B(u) = \tilde{s} \enspace \forall u \in [t,\Delta)
    \big] \\
    &=
    \int_{0}^\Delta {\mathbb E}\big[\tau^s_{k\ell}(d) \,|\, B(u) = s \enspace
    \forall u \in [0,t), B(u) = \tilde{s} \enspace \forall u \in [t,\Delta)
    \big]\frac{q_{ss'}e^{-q_{ss'}t}}{1 - e^{-q_{ss'}\Delta}}\,{\rm d}t \\ 
    &=
    \frac{q_{ss'}}{1 - e^{-q_{ss'}\Delta}} \int_{0}^\Delta (\Delta +  
    \phi_{s'}(d - v_{k\ell}(s)t - v_{k\ell}(s')(\Delta - t) \,|\,k,\ell))
    e^{-q_{ss'}t}\,{\rm d}t.
\end{align*}
Define $f(t,\Delta) := \phi_{s'}(d - v_{k\ell}(s)t - v_{k\ell}(s')(\Delta - t) \,|\,k,\ell)
    e^{-q_{ss'}t}$.
Now subtracting $\phi_{s}(d-v_{k\ell}(s)\Delta\,|\,k,\ell)$ from both sides,
dividing by $\Delta$ and letting $\Delta \downarrow 0$ we obtain
\begin{align*}
    v_{k\ell}(s)\phi_{s}'(d\,|\,k,\ell) &=
    1 + q_{ss}\phi_{s}(d\,|\,k,\ell) + \sum_{s' \neq s} q_{ss'} \lim_{\Delta \downarrow 0} \frac{q_{ss'}}{1 - e^{-q_{ss'}\Delta}} \int_{0}^\Delta f(t,\Delta)\,{\rm d}t \\
    &=
    1 + q_{ss}\phi_{s}(d\,|\,k,\ell) + \sum_{s' \neq s} q_{ss'} 
    \phi_{s'}(d\,|\,k,\ell),
\end{align*}
where the second step follows from L'Hopital's rule in combination with Leibniz' integral rule. 
We have thus obtained the desired system of linear differential equations:
\begin{align*}
    v_{k\ell}(s)\phi_{s}'(d\,|\,k,\ell) = 
    1 + \sum_{s' \in \mathcal{I}} q_{ss'}\phi_{s'}(d\,|\,k,\ell).
\end{align*}
The same steps can be performed to derive the second system of linear differential equations.
Again conditioning on a possible jump of the background process, as $\Delta\downarrow 0$,
\begin{align*}
    \psi_{ss'}(d\,|\,\alpha,k,\ell) 
    &= (1 + \Delta q_{ss})\,
    e^{-\alpha \Delta}\, \phi_{ss'}(d - v_{k\ell}(s)\Delta \,|\,k,\ell,\alpha) \\
    &+\sum_{\tilde{s} \neq s} \Delta q_{s\tilde{s}}
    \frac{q_{s\tilde{s}}e^{-\alpha\Delta}}{1 - e^{-q_{s\tilde{s}}\Delta}} \int_{0}^\Delta 
    \phi_{ss'}(d - v_{k\ell}(s)t - v_{k\ell}(\tilde{s})(\Delta - t) \,|\,k,\ell,\alpha)
    e^{-q_{s\tilde{s}}t}\,{\rm d}t + o(\Delta).  
\end{align*}
Subtracting $\psi_{ss'}(d-v_{k\ell}(s)\Delta\,|\,k,\ell,\alpha)$ from both sides and
expanding $e^{-\alpha\Delta}$ as $1-\alpha\Delta + o(\Delta)$ gives
\begin{align*}
    \psi_{ss'}(d\,|\,k,\ell,\alpha) - \psi_{ss'}(d-v_{k\ell}(s)\Delta\,|\,k,\ell,\alpha)
    &= (-\alpha\Delta + \Delta q_{ss} - \alpha\Delta^2 q_{ss}) \,\psi_{ss'}(d - v_{k\ell}(s)\Delta \,|\,k,\ell,\alpha) \\
    &\quad + \sum_{\tilde{s} \neq s} \frac{(1 - \alpha\Delta)\Delta q_{s\tilde{s}}^2}{1 - e^{-q_{s\tilde{s}}\Delta}} \int_{0}^\Delta f(t,\Delta)\, {\rm d}t + o(\Delta),
\end{align*}
with now $f(t,\Delta) := \phi_{ss'}(d - v_{k\ell}(s)t - v_{k\ell}(\tilde{s})(\Delta - t) \,|\,k,\ell,\alpha)
    e^{-q_{s\tilde{s}}t}$.
Dividing by $\Delta$, letting $\Delta\downarrow 0$ and using both
L'Hopital's and Leibniz' rule, we eventually obtain
\begin{align*}
    v_{k\ell}(s)\,\psi_{ss'}'(d\,|\,k,\ell,\alpha) &=
    (q_{ss} - \alpha)\,\psi_{ss'}(d\,|\,k,\ell,\alpha) + \sum_{\tilde{s} \neq s} q_{s\tilde{s}} \,
    \psi_{ss'}(d\,|\,k,\ell,\alpha).
\end{align*}
This completes the proof.
\end{proof}

\section{Upper bounds probabilities}
\label{appendix:c}
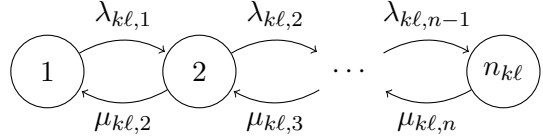
\begin{figure}[h]
\centering
\begin{tikzpicture}[shorten >=1pt,node distance=2cm,on grid,auto]
    \node[state] (q_0) {$1$};
    \node[state] (q_1) [right=of q_0] {$2$};
    \node[state, draw=none] (q_dots) [right=of q_1] {$\cdots$};
    \node[state] (q_4) [right=of q_dots] {$n_{k\ell}$};

    \path[->]
    (q_0) edge [bend left] node {$\lambda_{k\ell,1}$} (q_1)
    (q_1) edge [bend left] node {$\mu_{k\ell,2}$} (q_0)
    (q_1) edge [bend left] node {$\lambda_{k\ell,2}$} (q_dots)
    (q_dots) edge [bend left] node {$\mu_{k\ell,3}$} (q_1)
    (q_dots) edge [bend left] node {$\lambda_{k\ell,{n-1}}$} (q_4)
    (q_4) edge [bend left] node {$\mu_{k\ell,n}$} (q_dots);
\end{tikzpicture}
    \caption{Birth-death process}
    \label{fig:bd}
\end{figure}

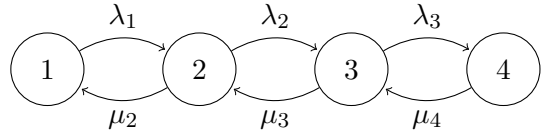
\begin{figure}[h]
    \centering
\begin{tikzpicture}[shorten >=1pt,node distance=2cm,on grid,auto]
    \node[state] (q_0) {$1$};
    \node[state] (q_1) [right=of q_0] {$2$};
    \node[state] (q_3) [right=of q_1] {$3$};
    \node[state] (q_4) [right=of q_3] {$4$};

    \path[->]
    (q_0) edge [bend left] node {$\lambda_{1}$} (q_1)
    (q_1) edge [bend left] node {$\mu_{2}$} (q_0)
    (q_1) edge [bend left] node {$\lambda_{2}$} (q_3)
    (q_3) edge [bend left] node {$\mu_{3}$} (q_1)
    (q_3) edge [bend left] node {$\lambda_{3}$} (q_4)
    (q_4) edge [bend left] node {$\mu_{4}$} (q_3);
\end{tikzpicture}
    \caption{Birth-death process with 4 states}
    \label{fig:bd4}
\end{figure}

We will show how results on the sum of independent exponentially random variables  can
be used to derive an upper bound on the hitting probability
\begin{align}
    \mathbb{P}(\exists t \in [0,M]:
    B_{k\ell}(t) = s_{k\ell}' \,|\, B_{k\ell}(0) = s_{k\ell}), \label{eq:ub2}
\end{align}
in case $B_{k\ell}(t)$ is a birth-death process.
Figure~\ref{fig:bd} shows the outline of such a process with $n_{k\ell}$ states.
By the structure of a birth-death process, it holds that
\begin{align*}
    \mathbb{P}(\exists t \in [0,M]:
    B_{k\ell}(t) = s_{k\ell}' \,|\, B_{k\ell}(0) = s_{k\ell}) \leqslant \mathbb{P}(S_{N} \leqslant M),
\end{align*}
with $S_{N}$ the sum of $N$ independent
exponentially distributed random variables.
Look for example at Figure~\ref{fig:bd4} with $\mathcal{S}_{k\ell} = \{1,2,3,4\}$
and note that 
\begin{align*}
    \mathbb{P}(\exists t \in [0,M]:
    B_{k\ell}(t) = 3 \,|\, B_{k\ell}(0) = 1) &\leqslant \mathbb{P}(E(\lambda_1) + E(\lambda_2) \leqslant M),
\end{align*}
with $E(\lambda_1) \sim \Exp(\lambda_1)$ and $E(\lambda_2) \sim \Exp(\lambda_2)$ independent.

With $S_N$
the sum of $N$ independent exponentially distributed random variables
with rates $\lambda_i$ 
(where $\lambda_i \neq \lambda_j$ for $i \neq j$) and density 
functions $f_i(z)=:\lambda_i\,e^{-\lambda_i z}$, we have \cite{Bibinger2013}
\begin{align*}
    f_{S_N}(z) &= \sum_{n = 1}^{N} \left(\prod_{j = 1,j \neq n}^N \frac{\lambda_j}{\lambda_j - 
    \lambda_n}\right) f_n(z) =: \sum_{n = 1}^N c_n f_n(z). 
\end{align*}
This implies for $M < \infty$ that
\begin{align*}
    \mathbb{P}(S_N \leqslant M) = \int_0^M \sum_{n = 1}^N c_n f_n(z)\,{\rm d}z = \sum_{n = 1}^N c_n
    \int_0^M f_n(z)\,{\rm d}z = \sum_{n = 1}^N c_n(1-e^{-\lambda_n M}).
\end{align*}
These results can now directly be used to find an upper bound on the probabilities in (\ref{eq:ub2}).

\end{document}